\documentclass{article}
\usepackage{amssymb}
\usepackage{mathrsfs}
\usepackage{cite}

\usepackage{stmaryrd}
\usepackage{amsfonts}
\usepackage{amsmath}
\usepackage{bm}
\usepackage{indentfirst}
\usepackage{array}
\usepackage{arydshln}
\usepackage{graphicx}
\usepackage{listings}
\usepackage{epstopdf}
\numberwithin{equation}{section}
\usepackage{geometry}
\usepackage{amsthm}
\newtheorem{theorem}{Theorem}

\newtheorem{lemma}{Lemma}
\newtheorem{rem}{Remark}
\newtheorem{proposition}{Proposition}
\title{\Large \bf{Large time behaviour of solutions to the $N-$dimensional scalar conservation law under periodic perturbations with nonlinear degenerate viscosity}
\setcounter{footnote}{-1}
\author{Yechi Liu$^\dag$
\noindent\footnote{\dag\quad E-mail: lyc9009@sina.cn.}\\
\small College of Science, National University of Defence Technology, Changsha 410003, P.R.China}}
\date{}
\begin{document}\large
\maketitle
\textbf{Abstract}. In this paper, we discuss the asymptotic behaviour of the weak solution to the Cauchy problem for the scalar viscous conservation law, with nonlinear Laplacian viscosity. Firstly, we obtain the existence, uniqueness and regularity of solutions when the initial data $u_0\in C^1(\mathbb R^N)\cap W^{1,\infty}(\mathbb R^N)$. Secondly, when $u_0$ is periodic, we prove the time-decay rate of the periodic solution and its gradient. At last, we study the long-time behaviour of perturbed solution to the Cauchy problem, in which the initial data is a $N-$d periodic perturbation around a planar rarefaction wave and obtain the time-decay rate of the perturbed solution approaching approximate planar rarefaction wave. The proof is given by technical energy methods and iteration technique.\\


\section{\Large Introduction and main results}
In this paper, we are concerned with a scalar conservation law with nonlinear Laplacian viscosity, which reads in $\mathbb R^N$ as
\begin{equation}\label{0}
\partial_tu(t,x)+\textrm{div}f\big(u(t,x)\big)=\textrm{div}\big(|\nabla u(t,x)|^{m-1}\nabla u(t,x)\big),\quad t\in(0,\infty), x\in\mathbb R^N,
\end{equation}
where $1m>1$ and the flux function $f=(f_1,\cdots,f_n)^T$ is smooth and satisfies $f(0)=f^\prime(0)=0,f_1^{\prime\prime}\geqslant c_f$ for some constant $c_f>0$.

In one-dimensional case with non-viscosity, i.e. $N=1$ and the right-hand side equals zero, the wave solutions to \eqref{0} include shocks and rarefaction waves (see \cite{lax57,liu78}). A centered rarefaction wave $u^R(t,x_1)$ is an entropy solution to the Riemann problem
\begin{equation}\label{R}
\left\{\begin{aligned}
&\partial_tu^R(t,x_1)+\partial_{x_1}f_1\big(u^R(t,x_1)\big)=0,\quad t\in(0,\infty),\quad x_1\in\mathbb R,\\
&u^R(0,x_1)=\left\{\begin{aligned}
    &u_-,\quad x_1<0,\\
    &u_+,\quad x_1>0,\end{aligned}\right.
\end{aligned}\right.
\end{equation}
where $u_-<u_+$ are constants, and it has an explicit formula as
\begin{equation}\label{R1}
u^R(t,x_1)=u^R\big(\frac{t}{x_1}\big)\equiv\left\{\begin{aligned}
&u_-, &&x_1<\lambda_-t,\\
&(\lambda)^{-1}\left(\frac{x_1}{t}\right), &&\lambda_-t\leqslant x_1\leqslant\lambda_+t,\\
&u_+, &&\lambda_+t<x_1,
\end{aligned}\right.
\end{equation}
where $\lambda=f_1^\prime,\lambda_\pm=\lambda(u_\pm)$. For the case with $m=1$, Il'in-Ole\u{\i}nik \cite{ili60} proved the Riemann solution consists of a single rarefaction wave solution, and the global solution in time tends toward the rarefaction wave. Hattori-Nishihara \cite{hat91} also showed the decay rate in time of the solution toward the single rarefaction wave in the $L^p$-norm ($p\in[1,\infty]$) for large $t>0$. Instead of initial data tending toward constants at far field states, Xin-Yuan-Yuan \cite{xin19} firstly investigated the large time behavior of the nonlinear waves under periodic perturbations for inviscid conservation laws, and they \cite{xin21} later extended the result to the linear viscosity case. When $m>1$, the nonlinear viscosity term in \eqref{0} models non-Newtonian fluid, such as blood, honey, butter, whipped cream, etc. (see \cite{lad70}). Such a viscosity term is also called the Ostwald-de Waele-type viscosity (see \cite{yos15,yos17} for more details). Because of the degeneracy of the viscosity, there were few results. Matsumura-Nishihara \cite{mat94} analysed asymptotic stability of a single rarefaction wave and Yoshida \cite{yos15} gave the time-decay rate. For the case with degenerate flux, Yoshida \cite{yos17a} found the asymptotic behaviour of the solution toward a multi-wave pattern including rarefaction waves and contact discontinuity as time tends to infinity, and obtained the time-decay rate in \cite{yos17}.

For the multiple dimensional case, Huang-Yuan \cite{hua21} proved that the solution to the scalar conservation law with linear viscosity time-asymptotically tends to the planar rarefaction waves when the initial perturbations are multi-dimensional periodic, and also obtained the time-decay rate. Recently, Huang-Xu-Yuan \cite{hua22} obtained asymptotic stability of planar rarefaction wave under 3-d periodic perturbations for Navier-Stokes equations. Concerning the case with degenerate viscosity, however, there is no result yet.

\vskip 0.2in

In this paper, for \eqref{0}, we consider the long-time behaviour of perturbed solution to the Cauchy problem, in which the initial data is a $N-$d periodic perturbation around a planar rarefaction wave. We want to use the periodic solutions and the planar rarefaction wave to construct an approximation to the planar rarefaction wave, and prove the perturbed solution converges to this approximation. For this purpose, we require the existence, regularity and time-decay rate of the periodic solutions given with periodic initial data. Considering the degeneracy of the viscosity in \eqref{0}, we need to begin with studying the existence and regularity of solutions to the Cauchy problem of equation \eqref{0} with large initial data
\begin{equation}\label{inda}
u(0,x)=u_0(x).
\end{equation}

A weak solution $u(t,x)$ to the Cauchy problem (\ref{0},\ref{inda}) is a measurable function defined in $[0,\infty)\times\mathbb R^N$ which satiesfies\\
\begin{equation*}
u(t,x)\in C\big(0,T;L^1(U)\big)\cap L^m\big(0,T;W^{1,m}(U)\big)
\end{equation*}
for some $T>0$ and any bounded open set $U\subset\mathbb R^N$, and
\begin{equation*}
\begin{aligned}
\int_{U}u(t,x)&\varphi(t,x){\rm d}x+\int_0^t\!\!\int_{U}f(u)\cdot\nabla\varphi{\rm d}x{\rm d}t\\
&=\int_{U}u_0(x)\varphi(0,x){\rm d}x+\int_0^t\!\!\int_{U}\Big(u\partial_t\varphi-\big(|\nabla u(t,x)|^{m-1}\nabla u(t,x)\big)\cdot\nabla\varphi\big){\rm d}x{\rm d}t
\end{aligned}
\end{equation*}
for all $t\in(0,T)$ and $\varphi(t,x)\in W^{1,\infty}\big(0,T;L^\infty(U)\big)\cap L^\infty\big(0,T;W_0^{1,\infty}(U)\big)$.

For the existence and regularity of solutions to the Cauchy problem (\ref{0},\ref{inda}), we obtain the following results.
\begin{theorem}[Existence]\label{ext}
If $m>1$ and $u_0\in C^1(\mathbb R^N)\cap W^{1,\infty}(\mathbb R^N)$, the Cauchy problem (\ref{0},\ref{inda}) admits a unique weak solution $u(t,x)$ such that, for any $T>0$,\\
(i), $\|u\|_{L^\infty([0,\infty)\times\mathbb R^N)}\leqslant\|u_0\|_{L^\infty(\mathbb R^N)}$,\\
(ii), $\nabla u$ is bounded on $[0,T]\times\mathbb R^N$,\\
(iii), $u$ is H\"{o}lder continuous on $[0,T]\times\mathbb R^N$ with index $\frac{1}{2},1$ and $\nabla u$ is locally H\"{o}lder continuous on $(0,\infty)\times\mathbb R^N$.
\end{theorem}

\vskip 0.2in

If $u_0$ is a periodic function, which means
\begin{equation}\label{up}
\left\{\begin{aligned}
&\partial_tu+\textrm{div}f(u)=\textrm{div}(|\nabla u|^{m-1}\nabla u),\\
&u(0,x)=\bar u+w_0(x),
\end{aligned}\right.
\end{equation}
where $w_0\in C^1(\mathbb R^N)$ is periodic with respect to each $x_i,i=1,2,\cdots,N$ on the $n$-dimensional torus $\mathbb T^N=:\prod_{i=1}^N[0,1]$ and satisfies
\begin{equation}\label{w}
\int_{\mathbb T^N}w_0(x)\textrm{d}x=0.
\end{equation}
By Theorem \ref{ext}, the solution $u$ to \eqref{up} is unique and satisfies (i)-(iii) in Theorem \ref{ext}. Besides,
\begin{theorem}\label{pu}
For $q\in[2,\infty]$, there exists a positive constant $C_q$, such that the periodic solution $u$ to \eqref{up} satisfies
\begin{equation}\label{pro4}
\|u(t,\cdot)-\bar u\|_{L^q(\mathbb T^N)}\leqslant C_q(1+t)^{-\frac{1}{m-1}},\quad \|\nabla u\|_{L^{m+1}(\mathbb T^N)}\leqslant C(1+t)^{-\frac{2}{(m-1)(m+1)}}
\end{equation}
for $q\in[2,\infty]$, and
\begin{equation}\label{pro5}
\|\nabla u\|_{L^q(\mathbb T^N)}\leqslant
\left\{\begin{aligned}
&C_q(1+t)^{-\frac{2\gamma_q}{(m-1)(m+1)}},\,\, m>2,N\in[2,4],\\
&C_q(1+t)^{-\frac{2}{(m-1)(m+1)}},\,\, m\in(1,2],N\geqslant2,
\end{aligned}\right.
\end{equation}
for $q\in(m+1,\infty]$, where $\gamma_q=\min\{1,\alpha_q\}$, $\alpha_q=\frac{(m+1)\big(2(q+1)+N(m-2)\big)}{(m+q-1)\big(N(m-2)-2m+2\big)}$.
\end{theorem}

\vskip 0.2in

Since the centered rarefaction wave given in \eqref{R1} is only Lipschitz continuous with respect of $x_1$, we need to construct a smooth viscous rarefaction wave $\tilde u^R(t,x_1)$, which is a solution to
\begin{equation}\label{tildeR}
\left\{\begin{aligned}
&\partial_t\tilde u^R(t,x_1)+\partial_1f_1\big(\tilde u^R(t,x_1)\big)=0,\quad t\in(0,\infty),\quad x_1\in\mathbb R,\\
&\tilde u^R(0,x_1)=\lambda^{-1}\left(\frac{\lambda_++\lambda_-}{2}+\frac{\lambda_+-\lambda_-}{2} \frac{e^{x_1}-e^{-x_1}}{e^{x_1}+e^{-x_1}}\right)=:\tilde u_0^R(x_1),\quad x_1\in\mathbb R
\end{aligned}\right.
\end{equation}
(see \cite{yos15}). It is easy to see that
\begin{equation*}
\lim_{x_1\rightarrow\pm\infty}\tilde u_0^R(x_1)=u_\pm.
\end{equation*}
In addition, we have the following lemma.
\begin{lemma}[\!\cite{yos15}\,]\label{yos}
The solution $\tilde u^R(t,x_1)$ to \eqref{tildeR} satisfies that, for any $t>0$,\\
(1). $u_-<\tilde u^R(t,x_1)<u_+$, $\lim_{t\rightarrow\infty}\big\|(\tilde u^R-u^R)(t,\cdot)\big\|_\infty=0$.\\
(2). $0<\partial_1\tilde u^R(t,x_1)\leqslant\min\left\{\frac{C}{t},u_+-u_-\right\}$.\\
(3). $\big\|\big(\tilde u^R(t,\cdot)-u_-\big)\big(\tilde u^R(t,\cdot)-u_+\big)\big\|_1\leqslant C(1+t)^\varepsilon$ for any $\varepsilon\in(0,1)$.\\
(4). $\big\|\partial_1\tilde u^R(t,\cdot)\big\|_p\leqslant C_p(1+t)^{-1+\frac{1}{p}}\leqslant C_p\min\left\{u_+-u_-,(u_+-u_-)^{\frac{1}{p}}(1+t)^{-1+\frac{1}{p}}\right\}$,

\quad\!$\big\|\partial_1^2\tilde u^R(t,\cdot)\big\|_p\leqslant C_p(1+t)^{-1}$\qquad for any $p\in[1,\infty]$.\\
(5). $\big|\tilde u^R(t,x_1)-u_+\big|\leqslant C_\varepsilon(1+t)^{-1+\varepsilon}e^{-\varepsilon|x_1-\lambda_+t|}$ for any $x_1\geqslant\lambda_+t$,

\quad\!$\big|\tilde u^R(t,x_1)-u_-\big|\leqslant C_\varepsilon(1+t)^{-1+\varepsilon}e^{-\varepsilon|x_1-\lambda_-t|}$ for any $x_1\leqslant\lambda_-t$,

\quad\!$\big|\tilde u^R(t,x_1)-u^R\left(\frac{x_1}{t}\right)\big|\leqslant C_\varepsilon(1+t)^{-1+\varepsilon}$ for any $\lambda_-t\leqslant x_1\leqslant\lambda_+t$.\\
(6). $\big\|\tilde u^R(t,\cdot)-u^R\left(\frac{\cdot}{t}\right)\big\|_p\leqslant C_{p,\varepsilon}(1+t)^{-1+\frac{1}{p}+\varepsilon}$ for any $0<\varepsilon<1$ and $p\in[1,\infty]$.\\
The $C_p$, $C_\varepsilon$ and $C_{p,\varepsilon}$ given above are positive constants depending on $p$ and $\varepsilon$.
\end{lemma}

Let $u_l(t,x)$ and $u_r(t,x)$ to be the solution of Cauchy problem \eqref{up} with initial data
\begin{equation*}
u_l(0,x)=u_-+w_0(x),\quad u_r(0,x)=u_++w_0(x),
\end{equation*}
respectively. Define
\begin{equation*}
w_l(t,x)=u_l(t,x)-u_-,\quad w_r(t,x)=u_r(t,x)-u_+.
\end{equation*}
Set
\begin{equation}\label{g}
g(t,x_1)=\frac{\tilde u^R(t,x_1)-u_-}{u_+-u_-},\quad t\in(0,\infty),\quad x_1\in\mathbb R.
\end{equation}
From Lemma \ref{yos}, we have
\begin{equation}\label{gest}
\begin{aligned}
&\|\partial_1g(t,\cdot)\|_{L^p(\mathbb R)}\leqslant C_p(1+t)^{-1+\frac{1}{p}},\qquad 1\leqslant p\leqslant\infty,\\
&\int_\mathbb Rg(t,x_1)\big(1-g(t,x_1)\big)\textrm{d}x_1\leqslant C_\varepsilon(1+t)^\varepsilon,\qquad 0<\varepsilon<1.
\end{aligned}
\end{equation}
Set
\begin{equation}\label{tildeu}
\tilde u=u_l(1-g)+u_rg.
\end{equation}
It is easy to see that $\tilde u$ is periodic with respect to $x_2,\cdots,x_n$ and
\begin{equation}\label{wlr}
\tilde u-\tilde u^R=w_l(1-g)+w_rg.
\end{equation}
Furthermore,
\begin{equation*}
\|\tilde u(t,\cdot)-\tilde u^R(t,\cdot)\|_\infty\leqslant C(1+t)^{-\frac{1}{m-1}},\quad \lim_{x_1\rightarrow\pm\infty}\tilde u(t,x)=u_\pm.
\end{equation*}

\vskip 0.2in

Noting that
\begin{equation*}
\big\|\tilde u^R(t,\cdot)-u^R\big(\frac{\cdot}{t}\big)\big\|_p\leqslant C(p,\varepsilon)(1+t)^{-1+\frac{1}{p}+\varepsilon}
\end{equation*}
for any $(\varepsilon,p)\in(0,1)\times[1,\infty]$(see (6) in Lemma \ref{yos}), it follows
\begin{equation*}
\|\tilde u-u^R\|_\infty\leqslant C(1+t)^{-\gamma},
\end{equation*}
where $\gamma=\min\{1-\varepsilon,\frac{1}{m-1}\}$ for any small $\varepsilon\in(0,1)$. Therefore, $\tilde u$ is an approximation to $u^R$ and $\tilde u^R$ in $\infty-$norm.

\vskip 0.2in

Consider the Cauchy problem
\begin{equation}\label{01}
\left\{\begin{aligned}
&\partial_tu(t,x)+\textrm{div}f\big(u(t,x)\big)=\textrm{div}\big(|\nabla u(t,x)|^{m-1}\nabla u(t,x)\big),\\
&u(0,x)=u_0(x):=\tilde u_0^R(x_1)+w_0(x).
\end{aligned}\right.
\end{equation}
Denote $\Omega=\mathbb R\times\mathbb T^{N-1}$. We have the following result.
\begin{theorem}[Time-decay rate]\label{tdr}
Let the periodic function $w_0\in C^1(\mathbb T^N)$, $1<m\leqslant\frac{3}{2}$. Then \eqref{01} has a unique weak solution $u(t,x)$ satisfying
\begin{equation}\label{res3}
\|u(t,\cdot)-\tilde u(t,\cdot)\|_{L^r(\Omega)}\leqslant C(1+t)^{-\frac{r-2}{r(3m+1)}}
\end{equation}
for any $r\geqslant2$, where $C>0$ depends only on $f$, $N$, $m$ and $\|u_0\|_{W^{1,\infty}(\Omega)}$.
\end{theorem}

\vskip 0.2in

Huang-Yuan \cite{hua21} established a Gagliardo-Nirenberg-type inequality on the domain $\Omega=\mathbb R\times\mathbb T^{N-1}(N\geqslant2)$ and used this inequality to obtain the time-decay rate in the $\infty-$norm, for the solution to \eqref{01} with $m=1$ asymptotically tending toward the planar viscous rarefaction wave, by estimating $\|\phi(t,\cdot)\|_p$ and $\|\nabla\phi(t,\cdot)\|_q$($p\in[1,\infty),q\in[2,\infty)$) of the perturbation $\phi=u-\tilde u$. The periodic solution to the Cauchy problem with periodic initial data has an exponential time-decay rate for the linear viscosity case(i.e., $m=1$ in \eqref{01}), but a polynomial time-decay rate for the case of $m>1$(see (\ref{pro4},\ref{pro5})). For the 1-d case, Yoshida \cite{yos15} gave the time-decay rate in the $p-$norm($p\in[2,\infty]$), for the solution to \eqref{01} with $N=1$, approaching to the rarefaction wave $\tilde u^R$, if the initial perturbation satisfies $w_0\in L^2(\mathbb R)$ and $\frac{{\rm d}}{{\rm d}x}w_0\in L^{m+1}(\mathbb R)$, the latter is used to obtain the $L^\infty-$estimate. Due to the degeneracy of the viscosity term($m>1$) and the periodic initial perturbation (large initial data, $w_0\notin L^2(\mathbb R^N),\nabla w_0\notin L^{m+1}(\mathbb R^N)$), we failed to get the $\infty-$norm estimate (similar to that in \cite{hua21}). Thus, we only obtain the time-decay rate in $L^r,r\geqslant2$ approaching the approximate planar rarefaction.

\vskip 0.2in

The rest of this paper will be organized as follows. In Section 2, we will give some useful inequalities and prove Theorem 1. Then, in Section 3, the asymptotic properties of periodic perturbation (i.e. the proof of Theorem 2) will be given. At last, we will obtain the time-decay rate (i.e. the proof of Theorem 3) in Section 4. The proof of Lemma 2 in Section 2 is placed in Appendix.

\vskip 0.2in

\noindent\textbf{Notations}. For derivatives, we use $\partial_t$ to represent $\frac{\partial}{\partial t}$ and $\partial_i$ to represent $\frac{\partial}{\partial x_i}$. For function spaces, $L^p=L^p(U)$ and $W^{1,p}=W^{1,p}(U)$ denote the usual Lebesgue space and $k-$th order Sobolev space on any region $U$ with norms $\|\cdot\|_p$ and $\|\cdot\|_{k,p}$, respectively, which means
\begin{equation*}
\|v\|_{p;U}=:\left(\int_U|v(x)|^p\,\textrm{d}x\right)^\frac{1}{p},\quad \|v\|_{1,p;U}=:\|v\|_{p;U}+\|\nabla v\|_{p;U},
\end{equation*}
where $\nabla v=(\partial_1v,\cdots,\partial_Nv)^T$. We omit the region in the notations of norms without misunderstanding. We also denote $\|\cdot\|=\|\cdot\|_2$ for simplicity. For constants, we use $c$ and $C$ to represent uncertain positive constants suitably small and large respectively. In particular, $c(a_1,a_2,\cdots)$ and $C(b_1,b_2,\cdots)$ represent that the constant $c$ and $C$ depend only on $a_1,a_2,\cdots$ and $b_1,b_2,\cdots$, respectively. We denote $\Omega=\mathbb R\times\mathbb T^{N-1}$ in the rest of this paper.

\section{\Large Existence and regularity of solutions}
In this section, we will prove Theorem 1. Because of the degeneracy, we need to regularize the viscosity and initial data. Firstly, we will state some used inequalities in the proof.
\begin{lemma}\label{abp}
For any $1\leqslant p\leqslant2$, $a,b\in\mathbb R^N$, it holds
\begin{equation}\label{ab2}
\big||a|^{p-1}a-|b|^{p-1}b\big|\leqslant(|a|^{p-1}+|b|^{p-1})|a-b|.
\end{equation}
For any $q\geqslant1$, $a,b\in\mathbb R^N$, it holds
\begin{equation}\label{ab1}
\big(|a|^{q-1}a-|b|^{q-1}b\big)\cdot(a-b)\geqslant c(q)\big(|a|^{q-1}+|b|^{q-1}+|a-b|^{q-1}\big)|a-b|^2.
\end{equation}
\end{lemma}
The proof of Lemma \ref{abp} will be given in Appendix. In the proof of existence, we also need the following imbedding theorem.
\begin{lemma}[P.62 in \cite{lad68}]\label{hpq}
There exists a constant $C=C(p,q,N)>0$ such that, for any function $v\in C^\infty\big(0,T;L^q(U)\big)\cap L^p\big(0,T;W_0^{1,p}(U)\big)$, where $U\subset\mathbb R^N$, it holds
\begin{equation*}
\int_0^T\!\!\int_U|v|^h{\rm d}x{\rm d}t\leqslant C\left\|\int_U|v|^q{\rm d}x\right\|_\infty^\frac{p}{N}\int_0^T\!\!\int_U|\nabla v|^p{\rm d}x{\rm d}t,
\end{equation*}
where $h=\frac{p(q+N)}{N}$.
\end{lemma}

Denote $M_1=\|u_0\|_\infty$, $M_2=\|\nabla u_0\|_\infty$. We can construct a sequence of smooth functions $\{u_{0n}\},n=1,2,\cdots$ which uniformly converges to $u_0$ and satisfies
\begin{equation*}
|u_{0n}|\leqslant M_1,\quad |\nabla u_{0n}|\leqslant M_2.
\end{equation*}
Consider the Cauchy problem
\begin{equation}\label{un}
\left\{\begin{aligned}
&\partial_tu_n+{\rm div}f(u_n)={\rm div}\left(\left(|\nabla u_n|^2+\frac{1}{n}\right)^\frac{m-1}{2}\nabla u_n\right),\\
&u_n(0,x)=u_{0n}(x).
\end{aligned}\right.
\end{equation}
It is easy to see that \eqref{un} has a unique classic solution satisfying
\begin{equation}\label{un0}
\|u_n\|_\infty\leqslant\|u_{0n}\|_\infty\leqslant M_1
\end{equation}
for each $n$ by using the Maximum principle(see \cite{lad68}). Next, we will prove
\begin{equation}\label{unp}
\|\nabla u_n\|_{L^\infty([0,T]\times B_{1,y})}\leqslant C(T)
\end{equation}
for any $T>0$ and $y\in\mathbb R$, where $B_{\rho,y}=\big\{x\in\mathbb R^N\big||x-y|<\rho\big\}$ and $C$ is independent of $n$.
\begin{proof}[Proof of \eqref{unp}]
Define
\begin{equation*}
v_n=|\nabla u_n|^2+\frac{1}{n},\qquad v_{0n}=|\nabla u_{0n}|^2+\frac{1}{n}.
\end{equation*}
Let $\xi_k(x)\in C^\infty(B_{y,1+\frac{1}{k}}),k=1,2,\cdots$ satisfy
\begin{equation*}
0\leqslant\xi_k\leqslant1,\,\, \xi_k=1\big(x\in B_{y,1+\frac{1}{k+1}}\big),\,\, \xi=0\big(x\in\partial B_{y,1+\frac{1}{k}}\big),\,\, |\nabla\xi_k|\leqslant2k(k+1).
\end{equation*}
Differentiating \eqref{un}$_1$ with respect to $x_j$ for any $j=1,\cdots,N$, multiplying the resultant equation with $\xi^2v_n^{\alpha_k}\partial_ju$ for $\alpha_k\geqslant0,k=1,2,\cdots$ to be determined below, and summing with respect to $j$ from $1$ to $N$, we have
\begin{equation*}
\partial_t(\xi^2v_n^{\alpha_k+1})=2(\alpha_k+1)\xi^2v_n^{\alpha_k}\partial_ju_n{\rm div}\Big(\partial_j(v_n^\frac{m-1}{2}\nabla u_n)-\partial_j\big(f(u_n)\big)\Big).
\end{equation*}
Integrating over $[0,T]\times B_{y,1+\frac{1}{k}}$ implies
\begin{equation}\label{vn1}
\begin{aligned}
&\frac{1}{2(\alpha_k+1)}\int_{B_{y,1+\frac{1}{k}}}\!\!\xi^2v_n^{\alpha_k+1}(T,x){\rm d}x+\int_0^T\!\!\int_{B_{y,1+\frac{1}{k}}}\!\!
\overbrace{\xi^2\sum_{i,j}\partial_i(v_n^{\alpha_k}\partial_ju_n)\partial_j(v_n^\frac{m-1}{2}\partial_iu_n)}^{J_{11}}{\rm d}x{\rm d}t\\
&=\frac{1}{2(\alpha_k+1)}\int_{B_{y,1+\frac{1}{k}}}\!\!\xi^2v_{n0}^{\alpha_k+1}{\rm d}x-2\int_0^T\!\!\int_{B_{y,1+\frac{1}{k}}}\!\!
\underbrace{\xi v_n^{\alpha_k}\sum_{i,j}\partial_i\xi\partial_j(v_n^\frac{m-1}{2}\partial_iu_n)\partial_ju_n}_{J_{12}}{\rm d}x{\rm d}t\\
&\qquad\qquad+\int_0^T\!\!\int_{B_{y,1+\frac{1}{k}}}\!\!\underbrace{\sum_{i,j}\partial_i(\xi^2v_n^{\alpha_k}\partial_ju_n)
\partial_j\big(f_i(u_n)\big)}_{J_{13}}{\rm d}x{\rm d}t.
\end{aligned}
\end{equation}
With direct calculation, it holds
\begin{equation*}
\begin{aligned}
&J_{11}=\xi^2\Big(\alpha_k\frac{m-1}{2}v_n^{\alpha_k+\frac{m-1}{2}-2}(\nabla v_n\cdot\nabla u_n)^2+\frac{1}{2}(\alpha_k+\frac{m-1}{2})v_n^{\alpha_k+\frac{m-1}{2}-1}|\nabla v_n|^2\\
&\hspace{3.3in}+v_n^{\alpha_k+\frac{m-1}{2}}\sum_{i,j}(\partial_{ij}u_n)^2\Big),\\
&J_{12}\leqslant C\xi|\nabla\xi|\Big(\frac{m-1}{2}v_n^{\alpha_k+\frac{m-1}{2}-1}|\nabla v_n\cdot\nabla u_n||\nabla u_n|+\frac{1}{2}v_n^{\alpha_k+\frac{m-1}{2}}|\nabla v_n|\Big),\\
&J_{13}\leqslant C(\alpha_k+1)\big(v^{\alpha_k+1}+(\alpha_k+1)v^{\alpha_k}|\nabla v|\big).
\end{aligned}
\end{equation*}
Here, the constant $C$ depends only on $f$ and $\|u_0\|_\infty$. Thus, from \eqref{vn1}, we can conclude
\begin{equation*}
\begin{aligned}
&\sup_{0\leqslant\tau\leqslant T}\int_{B_{y,1+\frac{1}{k}}}\!\!\xi^2v_n^{\alpha_k+1}(\tau,x){\rm d}x+\int_0^T\!\!\int_{B_{y,1+\frac{1}{k}}}\!\!\xi^2v_n^{\frac{2\alpha_k+m-3}{2}}|\nabla v_n|^2{\rm d}x{\rm d}t\\
&\leqslant C\int_0^T\!\!\int_{B_{y,1+\frac{1}{k}}}\!\!\Big(|\nabla\xi|^2v_n^{\frac{2\alpha_k+m+1}{2}}
+v_n^{\frac{2\alpha_k-m+3}{2}}\Big){\rm d}x{\rm d}t+4\int_{B_{y,1+\frac{1}{k}}}\!\!\xi^2v_{n0}^{\alpha_k+1}{\rm d}x.
\end{aligned}
\end{equation*}
by using Young's inequality. Set $w_n=v_n^{\frac{2\alpha_k+m+1}{4}}$ and $\lambda_k=\lambda(\alpha_k)=\frac{4(\alpha_k+1)}{2\alpha_k+m+1}$, and note that since $m>1$, $\frac{4}{m+1}<\lambda_k<2$ implies $\xi^\frac{4}{\lambda_k}\leqslant\xi^2$ for $k=1,2,\cdots$, then
\begin{equation}\label{vn2}
\begin{aligned}
&\sup_{0\leqslant\tau\leqslant T}\int_{B_{y,1+\frac{1}{k}}}\!\!(\xi^\frac{2}{\lambda_k}w_n)^{\lambda_k}(\tau,x){\rm d}x+\int_0^T\!\!\int_{B_{y,1+\frac{1}{k}}}\!\!(\xi^\frac{2}{\lambda_k}|\nabla w_n|)^2{\rm d}x{\rm d}t\\
&\qquad\qquad\leqslant Ck^2(k+1)^2\int_0^T\!\!\int_{B_{y,1+\frac{1}{k}}}\!\!(w_n^2+1){\rm d}x{\rm d}t+4\int_{B_{y,1+\frac{1}{k}}}\!\!\xi^2v_{n0}^{\alpha_k+1}{\rm d}x.
\end{aligned}
\end{equation}

On the other hand, using Lemma \ref{hpq} on $\xi^{\frac{2}{\lambda_k}}w_n$ by choosing $p=2,q=\lambda_k$, we have
\begin{equation}\label{vn5}
\begin{aligned}
&\int_0^T\!\!\int_{B_{y,1+\frac{1}{k+1}}}\!\!|\nabla u_n|^{s_{k+1}}{\rm d}x{\rm d}t<\int_0^T\!\!\int_{B_{y,1+\frac{1}{k}}}\!\!(\xi^\frac{2}{\lambda_k}w_n)^\frac{2(N+\lambda_k)}{N}{\rm d}x{\rm d}t\\
&\leqslant C\left(\sup_{0\leqslant t\leqslant T}\int_{B_{y,1+\frac{1}{k}}}\!\!(\xi^\frac{2}{\lambda_k}w_n)^{\lambda_k}{\rm d}x\right)^\frac{2}{N}\int_0^T\!\!\int_{B_{y,1+\frac{1}{k}}}\!\!\left|\nabla(\xi^\frac{2}{\lambda_k}w_n)\right|^2{\rm d}x{\rm d}t\\
&\leqslant C\left(\sup_{0\leqslant t\leqslant T}\int_{B_{y,1+\frac{1}{k}}}\!\!(\xi^\frac{2}{\lambda_k}w_n)^{\lambda_k}{\rm d}x+\int_0^T\int_{B_{y,1+\frac{1}{k}}}\!\!\big(\xi^\frac{2}{\lambda_k}|\nabla w_n|\big)^2{\rm d}x{\rm d}t\right)^{1+\frac{2}{N}}\\
&\hspace{2.3in}+Ck^2(k+1)^2\left(\int_0^T\!\!\int_{B_{y,1+\frac{1}{k}}}\!\!w_n^2{\rm d}x{\rm d}t\right)^{1+\frac{2}{N}},
\end{aligned}
\end{equation}
where we used Young's inequality and $s_{k+1}=m-1+2(1+\frac{2}{N})(\alpha_k+1)$. Comparing \eqref{vn2} and \eqref{vn5} implies
\begin{equation}\label{vnbu}
\begin{aligned}
&\|\nabla u_n\|_{s_{k+1};[0,T]\times B_{y,1+\frac{1}{k+1}}}^{s_{k+1}}\\
&\qquad\leqslant C\Big(k^2(k+1)^2\big(\|\nabla u_n\|_{s_k;[0,T]\times B_{y,1+\frac{1}{k}}}^{s_k}+T\big)\Big)^{1+\frac{2}{N}}+M_2^{2(\alpha_k+1)},
\end{aligned}
\end{equation}
Here, we used \eqref{un0} and $s_k=m-1+2(\alpha_k+1)$. Thus,
\begin{equation}\label{sk}
s_{k+1}=m-1+2(1+\frac{2}{N})^k(\alpha_1+1).
\end{equation}

\vskip 0.2in

Set $\alpha_1=0$, then we only need the estimate of $\|\nabla u_n\|_{m+1;[0,T]\times B_{y,2}}$. Let $\eta(x)\in C^\infty(B_{y,3})$ satisfy
\begin{equation*}
0\leqslant\eta\leqslant1,\quad \eta=1(x\in B_{y,2}),\quad \eta=0(x\in\partial B_{y,3}),\quad |\nabla\eta|\leqslant2.
\end{equation*}
Multiplying \eqref{un}$_1$ with $\eta^{m+1}u_n$ and integrating the resultant equation over $[0,T]\times B_{y,3}$, we have
\begin{equation}\label{vn3}
\begin{aligned}
&\int_{B_{y,3}}\!\!\eta^{m+1}u_n^2(T,x){\rm d}x+2\int_0^T\!\!\int_{B_{y,3}}\!\!v_n^\frac{m-1}{2}\nabla(\eta^{m+1}u_n)\cdot\nabla u_n{\rm d}x{\rm d}t\\
&\qquad\qquad=\int_{B_{y,3}}\!\!\eta^{m+1}u_{n0}^2{\rm d}x+2\int_0^T\!\!\int_{B_{y,3}}\!\!\eta^{m+1}u_nf^\prime(u_n)\cdot\nabla u_n{\rm d}x{\rm d}t.
\end{aligned}
\end{equation}
Note that
\begin{equation*}
v_n^\frac{m-1}{2}\nabla(\eta^{m+1}u_n)\cdot\nabla u_n=\eta^{m+1}v_n^\frac{m-1}{2}|\nabla u_n|^2+(m+1)\eta^mv_n^\frac{m-1}{2}u_n\nabla\eta\cdot\nabla u_n
\end{equation*}
and
\begin{equation*}
c(v_n^\frac{m+1}{2}-1)\leqslant v_n^\frac{m-1}{2}|\nabla u_n|^2,\quad |\nabla u_n|^{m+1}<v_n^\frac{m-1}{2}|\nabla u_n|^2,\quad v_n^\frac{m-1}{2}|\nabla u_n|\leqslant C(|\nabla u_n|^m+1).
\end{equation*}
We can conclude from \eqref{vn3} that
\begin{equation}\label{vn4}
\int_0^T\!\!\int_{B_{y,2}}\!\!v_n^\frac{m+1}{2}{\rm d}x{\rm d}t\leqslant C\int_0^T\!\!\int_{B_{y,3}}\!\!(|u_n|^{m+1}+1){\rm d}x{\rm d}t+4\int_{B_{y,3}}\!\!u_{n0}^2{\rm d}x,
\end{equation}
where we used Young's inequality and $C$ depends only on $f$ and $\|u_0\|_\infty$.

From \eqref{vnbu} and \eqref{vn4}, we can conclude that
\begin{equation*}
\|\nabla u_n\|_{s_{k+1};[0,T]\times B_{y,1+\frac{1}{k+1}}}^{s_{k+1}}\leqslant C\big(1+o(1)\big)\prod_{l=0}^{k-1}\left(2(k-l)^2(k-l+1)^2\right)^{(1+\frac{2}{N})^{l+1}}
\end{equation*}
as $k\rightarrow\infty$. Denote
\begin{equation*}
A_k=\left(\prod_{l=0}^{k-1}\left(2(k-l)^2(k-l+1)^2\right)^{(1+\frac{2}{N})^{i+1}}\right)^\frac{1}{s_{k+1}},
\end{equation*}
then
\begin{equation*}
\ln A_k=\frac{(1+\frac{2}{N})^k}{s_{k+1}}\sum_{l=0}^{k-1}(\frac{1}{1+\frac{2}{N}})^{k-l-1}
\big(2\ln(k-l)+2\ln(k-l+1)+\ln2\big).
\end{equation*}
It is easy to see that
\begin{equation*}
\lim_{k\rightarrow\infty}s_k=\infty,\quad \lim_{k\rightarrow\infty}\frac{(1+\frac{2}{N})^k}{s_{k+1}}=\frac{1}{2}
\end{equation*}
from \eqref{sk}, and
\begin{equation*}
\lim_{k\rightarrow\infty}\sum_{l=0}^{k-1}(\frac{1}{1+\frac{2}{N}})^{k-l-1}\ln(k-l+1)
<\int_1^\infty(\frac{1}{1+\frac{2}{N}})^{x-2}\ln(x+2)<\infty,
\end{equation*}
which implies $\lim_{k\rightarrow\infty}A_k<\infty$. Hence, \eqref{unp} holds true from the continuity of $\nabla u_n$.
\end{proof}

\vskip 0.2in

Since $y$ is arbitrarily given, we can conclude from \eqref{unp} that $u_n$ is Lipschitz continuous in $x$ on $[0,T)\times\mathbb R$ for any $T>0$, and, therefore, H\"{o}lder continuous in $t$ with index $\frac{1}{2}$, uniformly with respect to $n$. Thus, we may apply Arzela-Ascoli theorem, so that there exists a subsequence of $\{u_n\}$, still denoted as $\{u_n\}$, which uniformly converges to $u\in C^{\frac{1}{2},1}$ on any compact subset of $[0,\infty)\times\mathbb R^N$. In addition, from \eqref{un0}, we have
\begin{equation*}
\|u\|_\infty\leqslant M_1.
\end{equation*}
On the other hand, by using Theorem 1.1 in \cite{dib85}, $\nabla u_n$ is H\"{o}lder continuous, uniformly in $n$, on any compact subset of $(0,\infty)\times\mathbb R^N$. Furthermore, once by Arzela-Ascoli theorem, there exists a subsequence of $\{\nabla u_n\}$, still denoted as $\{\nabla u_n\}$, which converges to $\nabla u$, uniformly on any compact subset of $(0,\infty)\times\mathbb R^N$. Thus, $u$ is a weak solution of Cauchy problem (\ref{0},\ref{inda}) and $\nabla u$ is H\"{o}lder continuous on any compact subset of $(0,\infty)\times\mathbb R^N$ by Theorem 1.1 in \cite{dib85}. Furthermore, since $\|\nabla u_0\|\leqslant M_2$, the property of parabolic equation implies that, $\nabla u$ is bounded on $[0,T]\times\mathbb R^N$. Hence, $u$ is H\"{o}lder continuous on $[0,T]\times\mathbb R^N$ for any $T>0$.

\vskip 0.2in

It remains to check the uniqueness. Suppose $u_1,u_2$ are solutions with same initial data and set $v=u_1-u_2$. Then, we have
\begin{equation}\label{12}
\partial_tv+{\rm div}\big(f(u_1)-f(u_2)\big)={\rm div}\big(|\nabla u_1|^{m-1}\nabla u_1-|\nabla u_2|^{m-1}\nabla u_2\big).
\end{equation}
Define
\begin{equation*}
A_\gamma(x)=(1+|x|^2)^{-\gamma},
\end{equation*}
where $\gamma>0$ is a constant to be determined below, and let $\xi(x)$ be the truncation function given above in this section. Multiply \eqref{12} with $\xi^2A_\gamma v$, we have
\begin{equation}\label{121}
\begin{aligned}
&\frac{1}{2}\partial_t\big(\xi^2A_\gamma v^2\big)+\overbrace{\xi^2A_\gamma\nabla v\cdot\big(|\nabla u_1|^{m-1}\nabla u_1-|\nabla u_2|^{m-1}\nabla u_2\big)}^{J_{21}}+\sum_{i=1}^N\partial_i(\cdots)\\
&\left.\begin{aligned}
=&\sum_{i=1}^N\left(\int_0^v\big(f_i(w+u_2)-f_i(u_2)\big){\rm d}w-v\big(f_i(u_1)-f_i(u_2)\big)\right)\partial_i(\xi^2A_\gamma)\\
&\hspace{3cm}-\xi^2A_\gamma\sum_{i=1}^N\partial_iu_2\int_0^v\big(f_i^\prime(w+u_2)-f_i^\prime(u_2)\big){\rm d}w\\
&\hspace{3cm}-v\nabla(\xi^2A_\gamma)\cdot\big(|\nabla u_1|^{m-1}\nabla u_1-|\nabla u_2|^{m-1}\nabla u_2\big),
\end{aligned}\right\}J_{22}
\end{aligned}
\end{equation}
where
\begin{equation*}
\begin{aligned}
(\cdots)=&\xi^2A_\gamma v\sum_{i=1}^N\big(|\nabla u_1|^{m-1}\partial_iu_1-|\nabla u_2|^{m-1}\partial_iu_2\big)\\
&+\xi^2A_\gamma\sum_{i=1}^N\left(v\big(f_i(u_1)-f_i(u_2)\big)-\int_0^v\big(f_i(w+u_2)-f_i(u_2)\big){\rm d}w\right).
\end{aligned}
\end{equation*}
Motivated by \cite{dib89}, set
\begin{equation*}
\begin{aligned}
a_{ij}&(t,x)=\delta_{ij}\int_0^1\Big|\nabla\big(su_1+(1-s)u_2\big)\Big|^{m-1}{\rm d}s\\
&+(m-1)\int_0^1\Big|\nabla\big(su_1+(1-s)u_2\big)\Big|^{m-3}
\partial_i\big(su_1+(1-s)u_2\big)\partial_j\big(su_1+(1-s)u_2\big){\rm d}s
\end{aligned}
\end{equation*}
for $i,j=1,\cdots,N$, where $\delta_{ij}=\left\{\begin{aligned}1,&i=j,\\0,&i\neq j.\end{aligned}\right.$ Then it holds
\begin{equation*}
{\rm div}\big(|\nabla u_1|^{m-1}\nabla u_1-|\nabla u_2|^{m-1}\nabla u_2\big)=\sum_{i,j=1}^N\partial_i(a_{ij}\partial_jv)
\end{equation*}
and
\begin{equation*}
a_0|\chi|^2\leqslant\sum_{i,j=1}^Na_{ij}\chi_i\chi_j\leqslant ma_0|\chi|^2
\end{equation*}
for any $\chi=(\chi_1,\cdots,\chi_N)\in\mathbb R^N$, where
\begin{equation*}
a_0(t,x)=\int_0^1\Big|\nabla\big(su_1+(1-s)u_2\big)\Big|^{m-1}{\rm d}s.
\end{equation*}
Obviously, for any given $T>0$,
\begin{equation*}
0\leqslant a_0(t,x)\leqslant C(T),
\end{equation*}
when $t\in[0,T]$. Using $a_{ij}$ and $a_0$ given above, we can conclude that
\begin{equation}\label{j21}
J_{21}\geqslant c(m)\xi^2A_\gamma a_0|\nabla v|^2
\end{equation}
and
\begin{equation*}
\Big|\nabla(\xi^2A_\gamma)\cdot\big(|\nabla u_1|^{m-1}\nabla u_1-|\nabla u_2|^{m-1}\nabla u_2\big)\Big|\leqslant ma_0\big|\nabla(\xi^2A_\gamma)\big||\nabla v|.
\end{equation*}
On the other hand,
\begin{equation*}
\big|\nabla A_\gamma\big|\leqslant C(\gamma,N)A_\gamma.
\end{equation*}
Thus,
\begin{equation}\label{j22}
\big|J_{22}\big|\leqslant\varepsilon\xi^2A_\gamma a_0|\nabla v|^2+C\big(\xi^2+|\nabla\xi|\big)A_\gamma v^2,
\end{equation}
where $\varepsilon$ is a positive constant suitably small. Substituting \eqref{j21} and \eqref{j22} into \eqref{121}, integrating over $[0,t]\times\mathbb R^N$ for any $t\in[0,T]$, we have
\begin{equation*}
\int_{\mathbb R^N}A_\gamma(x)v(t,x)^2{\rm d}x\leqslant C\int_0^t\int_{\mathbb R^N}A_\gamma(x)v(\tau,x)^2{\rm d}x{\rm d}\tau
\end{equation*}
by choosing $\gamma>\frac{N}{2}$. Then Grownwall's inequality implies
\begin{equation*}
\int_{\mathbb R^N}A_\gamma(x)v(t,x)^2{\rm d}x=0,
\end{equation*}
since $v(0,x)=0$. Hence, $v=0$, which proves the uniqueness. Therefore, the proof of Theorem 1 is completed.

\vskip 0.2in

\begin{rem}
If, in addition, $\nabla u_0$ is H\"{o}lder continuous on $\mathbb R^N$, then similar to the analysis in \cite{lie93}, the H\"{o}lder continuity of $\nabla u$ may be claimed up to $t=0$.
\end{rem}

\section{\Large Properties of periodic solutions}
In this section, we will prove Theorem \ref{pu}. The region of integration with respect to $x$ is $\mathbb T^N$. Since the initial data given in\eqref{up} is periodic, the periodicity of solution is obvious by the uniqueness of weak solution given in Theorem \ref{ext}. Then, it remains to prove \eqref{pro4} and \eqref{pro5}.
\begin{proof}[Proof of \eqref{pro4} and \eqref{pro5}]
Multiplying \eqref{up}$_1$ by $u-\bar u$ and integrating the resultant equation with respect to $x$ over $\mathbb T^N$, we have
\begin{equation}\label{up1}
\frac{\textrm{d}}{\textrm{d}t}\|u-\bar u\|^2+2\|\nabla u\|_{m+1}^{m+1}=0.
\end{equation}
Since $u\in W^{1,m+1}(\mathbb T^N)$, the Poincar\'{e}'s inequality
\begin{equation}\label{upi}
c\|u-\bar u\|\leqslant \|\nabla u\|_{m+1}
\end{equation}
holds and implies
\begin{equation*}
\frac{\textrm{d}}{\textrm{d}t}\|u-\bar u\|^2+2c\|u-\bar u\|^{m+1}\leqslant0.
\end{equation*}
Thus,
\begin{equation}\label{up2}
\|u-\bar u\|\leqslant C(1+t)^{-\frac{1}{m-1}}
\end{equation}
and, therefore,
\begin{equation}\label{up3}
\int_0^t\|\nabla u\|_{m+1}^{m+1}\textrm{d}s\leqslant \|w_0\|^2
\end{equation}
for any $t>0$ by using \eqref{up1}.

\vskip 0.2in

Differentiating \eqref{up}$_1$ with respect to $x_j$ for $j=1,\cdots,n$, multiplying the resultant equations by $|\nabla u|^{m-1}\partial_ju$ respectively and summing from $1$ to $n$, we have
\begin{equation}\label{up4}
\begin{aligned}
\frac{1}{m+1}&\partial_t|\nabla u|^{m+1}+\nabla\big(\textrm{div}f(u)\big)\cdot\big(|\nabla u|^{m-1}\nabla u\big)\\
&=\textrm{div}\big(\textrm{div}(|\nabla u|^{m-1}\nabla u)|\nabla u|^{m-1}\nabla u\big)-\big(\textrm{div}(|\nabla u|^{m-1}\nabla u)\big)^2.
\end{aligned}
\end{equation}
With direct calculation, it holds
\begin{equation*}
\nabla\big(\textrm{div}f(u)\big)\cdot\big(|\nabla u|^{m-1}\nabla u\big)=\frac{m}{m+1}\sum_{i=1}^Nf_i^{\prime\prime}(u)|\nabla u|^{m+1}\partial_iu+\frac{1}{m+1}\sum_{i=1}^N\partial_i\big(f_i^\prime(u)|\nabla u|^{m+1}\big).
\end{equation*}
Then, integrating \eqref{up4} with respect to $x$ over $\mathbb T^n$ implies
\begin{equation}\label{up5}
\frac{\textrm{d}}{\textrm{d}t}\|\nabla u\|_{m+1}^{m+1}+(m+1)\int_{\mathbb T^N}\big(\textrm{div}(|\nabla u|^{m-1}\nabla u)\big)^2\textrm{d}x\leqslant C\|\nabla u\|_{m+1}^{m+1},
\end{equation}
where we used the fact that $f_i^{\prime\prime}(u)$ and $\partial_iu$ are bounded. Thus, integrating \eqref{up5} over $(0,t)$ and using \eqref{up3}, we can conclude
\begin{equation}\label{up6}
\|\nabla u\|_{m+1}^{m+1}+\int_0^t\!\!\int_{\mathbb T^N}\big(\textrm{div}(|\nabla u|^{m-1}\nabla u)\big)^2\textrm{d}x\textrm{d}s\leqslant C(\|\nabla w_0\|_{m+1}^{m+1}+\|w_0\|^2)<\infty.
\end{equation}
Here, we need $\nabla u_0=\nabla w_0$ and $\nabla w_0\in C^1(\mathbb T^N)$. On the other hand, comparing \eqref{up1} and \eqref{up5}, it holds
\begin{equation}\label{up7}
\frac{\textrm{d}}{\textrm{d}t}\big(\|u-\bar u\|^2+\|\nabla u\|_{m+1}^{m+1}\big)+\|\nabla u\|_{m+1}^{m+1}\leqslant0.
\end{equation}

\vskip 0.2in

If $\|\nabla u(t,\cdot)\|_{m+1}>1$ for any $t\in(0,t_0)$ with some $t_0>0$, then, from \eqref{up7}, we have
\begin{equation*}
\frac{\textrm{d}}{\textrm{d}t}\big(\|u-\bar u\|^2+\|\nabla u\|_{m+1}^{m+1}\big)+c\big(\|u-\bar u\|^2+\|\nabla u\|_{m+1}^{m+1}\big)\leqslant0,
\end{equation*}
where we used Poincar\'{e}'s inequality \eqref{upi} and the fact $\|\nabla u\|_{m+1}^2\leqslant\|\nabla u\|_{m+1}^{m+1}$ since $m+1>2$. Hence,
\begin{equation}\label{up8}
\|u-\bar u\|^2+\|\nabla u\|_{m+1}^{m+1}\leqslant\big(\|w_0\|^2+\|\nabla w_0\|_{m+1}^{m+1}\big)e^{-ct}\rightarrow0,\quad t\rightarrow\infty.
\end{equation}
\eqref{up8}, together with the continuity of $\|\nabla u(t,\cdot)\|_{m+1}$, implies that, there exists a constant $t_1\geqslant t_0$ such that $\|\nabla u(t_1,\cdot)\|_{m+1}=1$ and
\begin{equation*}
\|\nabla u(t,\cdot)\|_{m+1}<1,\quad t_1<t<t_1+\tau
\end{equation*}
for some small $\tau>0$. Then, on $(t_1,t_1+\tau)$ we have
\begin{equation*}
\frac{\textrm{d}}{\textrm{d}t}\big(\|u-\bar u\|^2+\|\nabla u\|_{m+1}^{m+1}\big)+c\big(\|u-\bar u\|^{m+1}+\|\nabla u\|_{m+1}^{(m+1)\frac{m+1}{2}}\big)\leqslant0
\end{equation*}
from \eqref{up7} by using Poincar\'{e}'s inequality \eqref{upi}, which means
\begin{equation}\label{up9}
\frac{\textrm{d}}{\textrm{d}t}\big(\|u-\bar u\|^2+\|\nabla u\|_{m+1}^{m+1}\big)+c\big(\|u-\bar u\|^2+\|\nabla u\|_{m+1}^{m+1}\big)^{\frac{m+1}{2}}\leqslant0.
\end{equation}
Integrating \eqref{up9} on $(t_1,t)$ implies
\begin{equation*}
\|u-\bar u\|^2+\|\nabla u\|_{m+1}^{m+1}\leqslant\left(\frac{2}{2\big(\|u(t_1,\cdot)-\bar u\|^2+1\big)+c(m-1)(t-t_1)}\right)^\frac{2}{m-1}\leqslant C(1+t)^{-\frac{2}{m-1}}.
\end{equation*}
Thus, for any $t>t_1$, it holds
\begin{equation}\label{upbu0}
\|u-\bar u\|\leqslant C(1+t)^{-\frac{1}{m-1}},\quad \|\nabla u\|_{m+1}\leqslant C(1+t)^{-\frac{2}{(m-1)(m+1)}}.
\end{equation}
If $\|\nabla u(t,\cdot)\|_{m+1}\leqslant1$ for any $t>0$, we can also obtain \eqref{upbu0} with similar discussion.

\vskip 0.2in

Next, we will estimate $\|u-\bar u\|_\infty$. For any $p\geqslant1$, denote $C_i,i=1,2,3$ to be positive constants independent of $p$, and by multiplying \eqref{up}$_1$ by $|u-\bar u|^{p-1}(u-\bar u)$ and integrating the resultant equation with respect to $x$ over $\mathbb T^n$, we have
\begin{equation*}
\frac{1}{p+1}\frac{\textrm{d}}{\textrm{d}t}\|u-\bar u\|_{p+1}^{p+1}+p\left(\frac{m+1}{p+m}\right)^{m+1}\left\|\nabla\big(|u-\bar u|^\frac{p-1}{m+1}(u-\bar u)\big)\right\|_{m+1}^{m+1}=0,
\end{equation*}
which implies
\begin{equation}\label{up10}
\begin{aligned}
\frac{\textrm{d}}{\textrm{d}t}\|u-\bar u\|_{p+1}^{p+1}+2C_1&\frac{p(p+1)}{(p+m)^{m+1}}\left\|\big(|u-\bar u|^\frac{p-1}{m+1}(u-\bar u)\big)\right\|_{1,m+1}^{m+1}\\
&\qquad=2C_1\frac{p(p+1)}{(p+m)^{m+1}}\|u-\bar u\|_{p+m}^{p+m}.
\end{aligned}
\end{equation}
Note that $|u-\bar u|^\frac{p-1}{m+1}(u-\bar u)$ may not vanish on the boundary and its mean value may not be equal to zero, we can not use Poincar\'{e}'s inequality directly. Set $v=|u-\bar u|^{\frac{p+m}{m+1}-1}(u-\bar u)=|u-\bar u|^{\frac{p-1}{m+1}}(u-\bar u)$, then
\begin{equation}\label{upbu}
\|u-\bar u\|_{p+m}^{p+m}=\|v\|_{m+1}^{m+1}.
\end{equation}
Using interpolation inequality, we have
\begin{equation}\label{up11}
\|v\|_{m+1}\leqslant C\|v\|_{1,m+1}^{\theta_1}\|v\|^{1-\theta_1},
\end{equation}
where $\theta_1=\frac{N(m-1)}{2(m+1)+N(m-1)}$. On the other hand, it holds
\begin{equation}\label{up12}
\|u-\bar u\|_{2\frac{p+m}{m+1}}\leqslant C\|u-\bar u\|^{\theta_2}\|u-\bar u\|_{p+1}^{1-\theta_2},
\end{equation}
in which $\theta_2=\frac{m-1}{p+m}$. Note that $\|v\|=\|u-\bar u\|_{2\frac{p+m}{m+1}}^\frac{p+m}{m+1}$, we can conclude from \eqref{upbu}, \eqref{up11} and \eqref{up12} that
\begin{equation}\label{up13}
\begin{aligned}
&\|u-\bar u\|_{p+m}^{p+m}\\
&\leqslant C\|v\|_{1,m+1}^{(m+1)\theta_1}\|u-\bar u\|^{(m+1)(1-\theta_1)\frac{p+m}{m+1}\theta_2}\|u-\bar u\|_{p+1}^{(m+1)(1-\theta_1)\frac{p+m}{m+1}(1-\theta_2)}\\
&=C\big\||u-\bar u|^\frac{p-1}{m+1}(u-\bar u)\big\|_{1,m+1}^{(m+1)\theta_1}\|u-\bar u\|^{(p+m)(1-\theta_1)\theta_2}\|u-\bar u\|_{p+1}^{(p+m)(1-\theta_1)(1-\theta_2)}\\
&\leqslant\varepsilon\big\||u-\bar u|^\frac{p-1}{m+1}(u-\bar u)\big\|_{1,m+1}^{m+1}+C\varepsilon^{-\frac{\theta_1}{1-\theta_1}}\|u-\bar u\|^{m-1}\|u-\bar u\|_{p+1}^{p+1}
\end{aligned}
\end{equation}
for any small $\varepsilon>0$. Substituting \eqref{up13} into \eqref{up10} and letting $\varepsilon=\frac{1}{2}$, we have
\begin{equation}\label{up15}
\begin{aligned}
\frac{\textrm{d}}{\textrm{d}t}\|u-\bar u\|_{p+1}^{p+1}+C_1&\frac{p(p+1)}{(p+m)^{m+1}}\big\||u-\bar u|^\frac{p-1}{m+1}(u-\bar u)\big\|_{1,m+1}^{m+1}\\
&\qquad\leqslant C_2\frac{p(p+1)}{(p+m)^{m+1}}\|u-\bar u\|^{m-1}\|u-\bar u\|_{p+1}^{p+1}.
\end{aligned}
\end{equation}
Set
\begin{equation*}
w(\tau,\cdot)=\big(u(\tau,\cdot)-\bar u\big)(1+t)^\frac{1}{m-1},\quad\tau=\ln(1+t).
\end{equation*}
Then, from \eqref{up15}, it holds
\begin{equation*}
\begin{aligned}
\frac{\textrm{d}}{\textrm{d}\tau}\|w\|_{p+1}^{p+1}+ C_1&\frac{p(p+1)}{(p+m)^{m+1}}\big\||w|^\frac{p-1}{m+1}w\big\|_{1,m+1}^{m+1}\\
&\qquad\leqslant\frac{p+1}{m-1}\|w\|_{p+1}^{p+1}+C_2\frac{p(p+1)}{(p+m)^{m+1}}\|w\|^{m-1}\|w\|_{p+1}^{p+1}.
\end{aligned}
\end{equation*}
Since $\|w(\tau,\cdot)\|\leqslant C$ for any $\tau>0$ from \eqref{upbu0} and $\|w(0,\cdot)\|_\infty=\|u_0-\bar u\|_\infty<\infty$, we can conclude that
\begin{equation}\label{up16}
\frac{\textrm{d}}{\textrm{d}\tau}\|w\|_{p+1}^{p+1}+ C_1\frac{p(p+1)}{(p+m)^{m+1}}\big\||w|^\frac{p-1}{m+1}w\big\|_{1,m+1}^{m+1}\leqslant C_3(p+1)\|w\|_{p+1}^{p+1}
\end{equation}
for any $p\geqslant1$, where $0<C_1<1$ and $C_3>1$.

\vskip 0.2in

Set
\begin{equation*}
\begin{aligned}
\lambda_0=0,\quad&\lambda_k=(m+1)\lambda_{k-1}+1,\quad a_k=C_1\frac{\lambda_k(\lambda_k+1)}{(\lambda_k+m)^{m+1}},\quad b_k=C_3(\lambda_k+1),\\
&\alpha_k=\frac{\lambda_k+1}{\lambda_{k-1}+1},\quad\beta_k=\frac{(m+1)(N+\alpha_k)-N\alpha_k}{(m+1)(N+1)-N\alpha_k},\quad k=1,2,\cdots.
\end{aligned}
\end{equation*}
It is easy to see that, for any $k\geqslant2$,
\begin{equation*}
\lambda_k>1,\quad 0<a_k<1,\quad b_k>1,\quad 1<\alpha_k<m+1,\quad 1<\beta_k<\alpha_k.
\end{equation*}
Choose $p=\lambda_k$ and let $v_k=|w|^{\lambda_{k-1}+1}$ in \eqref{up16}, then
\begin{equation}\label{up17}
\frac{\textrm{d}}{\textrm{d}\tau}\|v\|_{\alpha_k}^{\alpha_k}+a_k\|v\|_{1,m+1}^{m+1}\leqslant b_k\|v\|_{\alpha_k}^{\alpha_k}.
\end{equation}
Interpolation inequality and Young's inequality imply that
\begin{equation}\label{up18}
\|v\|_{\alpha_k}^{\alpha_k}\leqslant C\|v\|_{1,m+1}^{\alpha_k\tilde\theta_k}\|v\|_1^{\alpha_k(1-\tilde\theta_k)}\leqslant \varepsilon_k\|v\|_{1,m+1}^{m+1}+C\varepsilon_k^{-\gamma_k}\|v\|_1^{\beta_k},
\end{equation}
where
\begin{equation*}
\varepsilon_k=\frac{a_k}{2b_k},\quad\tilde\theta_k=\frac{N(m+1)(\alpha_k-1)}{(Nm+m+1)\alpha_k}\in(0,1),\quad \gamma_k=\frac{\beta_k\tilde\theta_k}{(m+1)(1-\tilde\theta_k)}.
\end{equation*}
Note that $\varepsilon_k(b_k+\varepsilon_k)<a_k$, then multiplying \eqref{up18} by $b_k+\varepsilon_k$ and adding \eqref{up17} imply that
\begin{equation*}
\frac{\textrm{d}}{\textrm{d}\tau}\|v\|_{\alpha_k}^{\alpha_k}+ \big(a_k-\varepsilon_k(b_k+\varepsilon_k)\big)\|v\|_{1,m+1}^{m+1}+\varepsilon_k\|v\|_{\alpha_k}^{\alpha_k}\leqslant C(b_k+\varepsilon_k)\varepsilon_k^{-\gamma_k}\|v\|_1^{\beta_k},
\end{equation*}
which means
\begin{equation*}
\frac{\textrm{d}}{\textrm{d}\tau}\|w\|_{\lambda_k+1}^{\lambda_k+1}+\varepsilon_k\|w\|_{\lambda_k+1}^{\lambda_k+1}\leqslant \tilde C_k(b_k+\varepsilon_k)\left(\sup_{\tau\geqslant0}\|w\|_{\lambda_{k-1}+1}^{\lambda_{k-1}+1}\right)^{\beta_k}.
\end{equation*}
Note that $(\lambda_{k-1}+1)\beta_k<\lambda_k+1$ and when $k\rightarrow\infty$,
\begin{equation*}
\varepsilon_k\sim\frac{C}{(\lambda_k+1)^m},\quad b_k\sim C(\lambda_k+1),\quad \tilde C_k\sim C\varepsilon_k^{-\gamma_k}<C(\lambda_k+m)^\frac{Nm^2}{m+1}.
\end{equation*}
Thus, by using Lemma 3.2 in \cite{ali79}, we can conclude that
\begin{equation*}
\|w(\tau,\cdot)\|_\infty\leqslant C\max\left\{1,\sup_{\tau\geqslant0}\|w(\tau,\cdot)\|,\|w(0,\cdot)\|_\infty\right\},
\end{equation*}
which implies immediately
\begin{equation}\label{up19}
\|u(t,\cdot)-\bar u\|_\infty\leqslant C(1+t)^{-\frac{1}{m-1}}.
\end{equation}

\vskip 0.2in

Next, we will estimate $\|\nabla u\|_q$ for any $q>m+1$. Multiplying \eqref{up}$_1$ with $\textrm{div}\big(|\nabla u|^{q-2}\nabla u\big)$ and integrating with respect to $x$ over $\mathbb T^N$, we have
\begin{equation}\label{up20}
\begin{aligned}
\frac{1}{q}\frac{\textrm{d}}{\textrm{d}t}\|\nabla u\|_q^q+&\int_{\mathbb T^N}\textrm{div}\big(|\nabla u|^{m-1}\nabla u\big)\textrm{div}\big(|\nabla u|^{q-2}\nabla u\big)\textrm{d}x\\
&\qquad\qquad\qquad=\int_{\mathbb T^n}\textrm{div}\big(f(u)\big)\textrm{div}\big(|\nabla u|^{q-2}\nabla u\big)\textrm{d}x.
\end{aligned}
\end{equation}
Since for any $j=1,\cdots,n$,
\begin{equation*}
\begin{aligned}
&f_j^\prime(u)\partial_ju\textrm{div}\big(|\nabla u|^{q-2}\nabla u\big)\\
&=\sum_{i=1}^N\Big(\partial_i\big(f_j^\prime(u)\partial_ju|\nabla u|^{q-2}\partial_iu\big)-f_j^{\prime\prime}(u)|\nabla u|^{q-2}(\partial_iu)^2\partial_ju-f_j^\prime(u)|\nabla u|^{q-2}\partial_iu\partial_{ij}u\Big)
\end{aligned}
\end{equation*}
and
\begin{equation*}
q\sum_{i=1}^Nf_j^\prime(u)|\nabla u|^{q-2}\partial_iu\partial_{ij}u=f_j^\prime\partial_j\big(|\nabla u|^q\big)=\partial_j\big(f_j^\prime(u)|\nabla u|^q\big)-f_j^{\prime\prime}(u)|\nabla u|^q\partial_ju,
\end{equation*}
it holds
\begin{equation}\label{up21}
\begin{aligned}
&\int_{\mathbb T^N}\textrm{div}\big(f(u)\big)\textrm{div}\big(|\nabla u|^{q-2}\nabla u\big)\textrm{d}x\\
&\qquad\qquad=\frac{1}{q}\int_{\mathbb T^N}\sum_{j=1}^Nf_j^{\prime\prime}(u)|\nabla u|^q\partial_ju\textrm{d}x-\int_{\mathbb T^N}\sum_{i,j=1}^Nf_j^{\prime\prime}(u)|\nabla u|^{q-2}(\partial_iu)^2\partial_ju.
\end{aligned}
\end{equation}
In addition, we can conclude from \cite{nak00} that
\begin{equation}\label{up22}
\begin{aligned}
&\int_{\mathbb T^N}\textrm{div}\big(|\nabla u|^{m-1}\nabla u\big)\textrm{div}\big(|\nabla u|^{q-2}\nabla u\big)\textrm{d}x\geqslant\int_{\mathbb T^N}|\nabla u|^{m+q-3}|D^2u|^2\textrm{d}x\\
&\quad+\frac{m+q-5}{4}\int_{\mathbb T^N}|\nabla u|^{m+q-5}\big|\nabla(|\nabla u|^2)\big|^2\textrm{d}x-C\int_{\partial\mathbb T^N}|\nabla u|^{m+q-1}\textrm{d}S,
\end{aligned}
\end{equation}
where $|D^2u|^2=\sum_{i,j=1}^N(\partial_{ij}u)^2$. Furthermore, the trace theorem, the interpolation inequality and the Cauchy's inequality imply that
\begin{equation}\label{up23}
\begin{aligned}
C\int_{\partial\mathbb T^N}&|\nabla u|^{m+q-1}\textrm{d}S\leqslant C\big\||\nabla u|^\frac{m+q-1}{2}\big\|_{H^\frac{1}{2}}^2\\
&\leqslant C\big\||\nabla u|^\frac{m+q-1}{2}\big\|_{1,2}\big\||\nabla u|^\frac{m+q-1}{2}\big\|\leqslant\frac{1}{2}\big\||\nabla u|^\frac{m+q-1}{2}\big\|_{1,2}^2+C\|\nabla u\|_{m+q-1}^{m+q-1}.
\end{aligned}
\end{equation}
Substituting \eqref{up21}-\eqref{up23} into \eqref{up20}, it holds
\begin{equation}\label{up24}
\frac{\textrm{d}}{\textrm{d}t}\|\nabla u\|_q^q+C_4\big\||\nabla u|^\frac{m+q-1}{2}\big\|_{1,2}^2\leqslant C_5q^2\|\nabla u\|_{m+q-1}^{m+q-1}+C_6q\|\nabla u\|_{q+1}^{q+1},
\end{equation}
where $C_i,i=4,5,6$ are positive constants independent of $q$.

\vskip 0.2in

\noindent\textbf{Case 1}:$m\leqslant2$ and $N\in[2,4]$.

We can use Lemma 1 in \cite{nak00} to conclude that
\begin{equation}\label{up25}
\|\nabla u\|_{m+q-1}^{m+q-1}\leqslant\varepsilon\big\||\nabla u|^\frac{m+q-1}{2}\big\|_{1,2}^2+C\|\nabla u\|_{m+1}^{m+q-1}
\end{equation}
and
\begin{equation}\label{up251}
\|\nabla u\|_{q+1}^{q+1}\leqslant\varepsilon\big\||\nabla u|^\frac{m+q-1}{2}\big\|_{1,2}^2+C\|\nabla u\|_{m+1}^{\alpha_q(m+q-1)}
\end{equation}
for some small $\varepsilon>0$, where
\begin{equation*}
\alpha_q=\frac{(1-\theta_3)(q+1)}{(1-\theta_3)(q+1)+m-2}\leqslant0,\quad \theta_3=\frac{m+q-1}{2}\cdot\frac{\frac{1}{m+1}-\frac{1}{q+1}}{\frac{1}{n}-\frac{1}{2}+\frac{m+q-1}{2(m+1)}},
\end{equation*}
so that $\alpha_q=\frac{(m+1)\big(2(q+1)+N(m-2)\big)}{(m+q-1)\big(N(m-2)-2m+2\big)}$. It follows from \eqref{up24} and \eqref{up25} that
\begin{equation}\label{up26}
\frac{\textrm{d}}{\textrm{d}t}\|\nabla u\|_q^q+\big\||\nabla u|^\frac{m+q-1}{2}\big\|_{1,2}^2\leqslant C(q)\left(\|\nabla u\|_{m+1}^{m+q-1}+\|\nabla u\|_{m+1}^{\alpha_q(m+q-1)}\right).
\end{equation}
Again, by using Lemma 1 in \cite{nak00}, we have
\begin{equation*}
\|\nabla u\|_q^{m+q-1}\leqslant\varepsilon\big\||\nabla u|^\frac{m+q-1}{2}\big\|_{1,2}^2+C\|\nabla u\|_{m+1}^{m+q-1},
\end{equation*}
which, together with \eqref{up26}, implies
\begin{equation}\label{up261}
\begin{aligned}
\frac{\textrm{d}}{\textrm{d}t}\|\nabla u\|_q^q+c\|\nabla u\|_q^{m+q-1}\leqslant&C\left(\|\nabla u\|_{m+1}^{m+q-1}+\|\nabla u\|_{m+1}^{\alpha_q(m+q-1)}\right)\\
\leqslant&C\left((1+t)^{-\frac{2(m+q-1)}{(m-1)(m+1)}}+(1+t)^{-\frac{2\alpha_q(m+q-1)}{(m-1)(m+1)}}\right),
\end{aligned}
\end{equation}
where we used \eqref{upbu0}. Thus, from Lemma 3 in \cite{nak00}, it holds
\begin{equation}\label{up27}
\|\nabla u\|_q\leqslant C(1+t)^{-\frac{2\gamma_q}{(m-1)(m+1)}}
\end{equation}
for any $q\geqslant2$, where $\gamma_q=\min\{1,\alpha_q\}$.

\vskip 0.2in

\noindent\textbf{Case 2}:$m\in(1,2)$ and $N\geqslant2$.

Since there are some $1<m<2$ making $\alpha_q<0$ when the space dimension $N>4$, we can not obtain the estimate \eqref{up251} by Young's inequality. However, since $\|\nabla u\|_\infty$ is bounded and $m<2$, we have
\begin{equation*}
\|\nabla u\|_{q+1}^{q+1}\leqslant C\|\nabla u\|_{m+q-1}^{m+q-1}
\end{equation*}
instead of \eqref{up251}. Thus,
\begin{equation*}
\frac{\textrm{d}}{\textrm{d}t}\|\nabla u\|_q^q+c\|\nabla u\|_q^{m+q-1}\leqslant C\|\nabla u\|_{m+1}^{m+q-1}\leqslant C(1+t)^{-\frac{2(m+q-1)}{(m-1)(m+1)}},
\end{equation*}
which implies
\begin{equation}\label{up28}
\|\nabla u\|_q\leqslant C(1+t)^{-\frac{2}{(m-1)(m+1)}}.
\end{equation}
Since $\nabla u(t,\cdot)$ is continuous on $\mathbb T^N$, \eqref{up27} and \eqref{up28} holds true for $q=\infty$.
\end{proof}

\section{\Large Time-decay rate}
In this section, we will prove Theorem \ref{tdr}. Firstly, we need the following lemma.
\begin{lemma}\label{l1est}
Let $u,v$ be the solutions to \eqref{0} with initial data $u_0,v_0$ respectively, where $u_0$ and $v_0$ are periodic in $x_i,i=2,\cdots,N$, satisfy the condition in Theorem \ref{ext} and, in addition, $u_0-v_0\in L^1(\Omega)$. Then
\begin{equation}\label{uv0}
\|u-v\|_1\leqslant\|u_0-v_0\|_1.
\end{equation}
\end{lemma}
\begin{proof}
Set
\begin{equation}\label{jdel}
J_\delta(\eta)=(\eta^2+\delta^2)^\frac{1}{2},
\end{equation}
where $\delta>0$ is any constant, then
\begin{equation*}
J_\delta^\prime(\eta)=\frac{\eta}{(\eta^2+\delta^2)^{\frac{1}{2}}},\quad J_\delta^{\prime\prime}(\eta)=\frac{\delta^2}{(\eta^2+\delta^2)^{\frac{3}{2}}}.
\end{equation*}
Let $\xi_n(x_1)$ be a truncation function which equals $1$ on $[-n,n]$ and vanishes out of $(-n-1,n+1)$ and satisfies $|\partial_1\xi_n|\leqslant2$. Since $u,v$ are solutions, we have
\begin{equation}\label{uv1}
\partial_t(u-v)+{\rm div}\big(f(u)-f(v)\big)={\rm div}\big(|\nabla u|^{m-1}\nabla u-|\nabla v|^{m-1}\nabla v\big).
\end{equation}
Multiplying \eqref{uv1} with $\xi_nJ_\delta^\prime(u-v)$ and integrating with respect to $x$ over $\Omega$ imply
\begin{equation}\label{uv2}
\begin{aligned}
&\partial_t\int_\Omega\xi_nJ_\delta(u-v){\rm d}x+\overbrace{\int_\Omega\xi_nJ_\delta^{\prime\prime}(u-v)\nabla(u-v)\cdot \big(|\nabla u|^{m-1}\nabla u-|\nabla v|^{m-1}\nabla v\big){\rm d}x}^{I_{01}}\\
&=\underbrace{\int_\Omega\xi_nJ_\delta^{\prime\prime}(u-v)\nabla(u-v)\cdot\big(f(u)-f(v)\big){\rm d}x}_{I_{02}}+\underbrace{\int_\Omega J_\delta^\prime(u-v)\partial_1\xi_n\cdot K(u,v){\rm d}x}_{I_{03}},
\end{aligned}
\end{equation}
where
\begin{equation*}
K(u,v)=\big(f(u)-f(v)\big)-\big(|\nabla u|^{m-1}\nabla u-|\nabla v|^{m-1}\nabla v\big).
\end{equation*}
It is easy to see that, for any $T>0$,
\begin{equation*}
\int_0^TI_{01}{\rm d}t\geqslant c(m)\int_0^T\!\!\int_\Omega\xi_nJ_\delta^{\prime\prime}(u-v)|\nabla(u-v)|^{m+1}{\rm d}x{\rm d}t\geqslant 0
\end{equation*}
from \eqref{ab1}, and
\begin{equation*}
\lim_{\delta\rightarrow0}\int_0^TI_{02}{\rm d}t=0.
\end{equation*}
Since $u$, $v$, $\nabla u$ and $\nabla v$ are bounded on $[0,T]\times\Omega$, we have
\begin{equation*}
\int_0^TI_{03}{\rm d}t\leqslant C(T)\int_{([-n-1,-n]\cup[n,n+1])\times\mathbb T^{N-1}}|\partial_1\xi_n|{\rm d}x\leqslant C(T).
\end{equation*}
Thus, integrating \eqref{uv2} over $(0,T)$ and let $\delta\rightarrow0$ imply
\begin{equation*}
\int_\Omega\xi_n|u-v|{\rm d}x\leqslant C(T)+\int_\Omega\xi_n|u_0-v_0|{\rm d}x.
\end{equation*}
Therefore, $u-v\in L^1(\Omega)$ for any $T>0$ by letting $n\rightarrow\infty$. On the other hand, since $u$ and $v$ are H\"{o}lder continuous on $[0,T]\times\Omega$, it holds
\begin{equation*}
\lim_{x_1\rightarrow\pm\infty}(u-v)=0.
\end{equation*}
Then, multiply \eqref{uv1} with $J_\delta^\prime(u-v)$ and integrate with respect to $x$ over $\Omega$. From the discussion about $I_{01}$ and $I_{02}$ above, we can conclude that
\begin{equation*}
\int_\Omega|u-v|{\rm d}x\leqslant\int_\Omega|u_0-v_0|{\rm d}x,
\end{equation*}
which completes the proof.
\end{proof}

\begin{rem}\label{hal}
Similar to the discussion above, we can conclude that, if $u_0-v_0\in L^1(\mathbb R_+\times\mathbb T^{N-1})$ instead of $u_0-v_0\in L^1(\Omega)$ in Lemma \ref{l1est}, then $u-v\in L^1(\mathbb R_+\times\mathbb T^{N-1})$ by choosing $\xi_n(x)$ equals 1 on $[0,n]$ and vanishes on $(n+1,\infty)$. The same conclusion also holds true with $\mathbb R_+\times\mathbb T^{N-1}$ replaced by $\mathbb R_-\times\mathbb T^{N-1}$.
\end{rem}

\vskip 0.2in

In the rest of this section, let $1<m\leqslant\frac{3}{2}$. Define $\phi(t,x)=u(t,x)-\tilde u(t,x)$, then $\phi$ is periodic with respect to $x_2,\cdots,x_N$ with period 1 and $\phi(0,x)=0$. From \eqref{tildeR}$_1$, \eqref{g}, \eqref{tildeu} and \eqref{01}$_1$, we can conclude that
\begin{equation}\label{phi}
\begin{aligned}
&\partial_t\phi-{\rm div}\big(|\nabla(\phi+\tilde u)|^{m-1}|\nabla(\phi+\tilde u)-|\nabla\tilde u|^{m-1}\nabla\tilde u\big)\\
&\qquad\qquad=-{\rm div}\big(f(\phi+\tilde u)-f(\tilde u)\big)+{\rm div}(|\nabla\tilde u|^{m-1}\nabla\tilde u)-h,
\end{aligned}
\end{equation}
where
\begin{equation}\label{h}
h=\partial_t\tilde u+\textrm{div}f(\tilde u).
\end{equation}
In addition, since
\begin{equation*}
\int_\Omega|\phi|{\rm d}x\leqslant C\left(\int_{\mathbb R_+\times\mathbb T^{N-1}}+\int_{\mathbb R_-\times\mathbb T^{N-1}}\right)\big|(1-g)(u-u_l)+g(u-u_r)\big|{\rm d}x,
\end{equation*}
we can conclude that $\phi(t,\cdot)\in L^1(\Omega)$ for any $t>0$ from Remark \ref{hal}. Therefore,
\begin{equation*}
\lim_{x_1\rightarrow\pm\infty}\phi(t,x)=0
\end{equation*}
for any $t>0$.

With direct calculation, we have
\begin{equation*}
\begin{aligned}
h=&\overbrace{(1-g)g(u_r-u_l)\sum_{i=1}^N\big(\sigma_i(\tilde u,u_l)\partial_iw_l-\sigma_i(\tilde u,u_r)\partial_iw_r\big)+(u_r-u_l)(\tilde u-u^R)\sigma_1(\tilde u,u^R)\partial_1g}^{J_1}\\
  &\hspace{2cm}+\underbrace{(1-g)\textrm{div}(|\nabla u_l|^{m-1}\nabla u_l)+g\textrm{div}(|\nabla u_r|^{m-1}\nabla u_r)}_{J_2}
\end{aligned}
\end{equation*}
where
\begin{equation*}
\sigma_i(u,v)=\int_0^1f_i^{\prime\prime}\big(v+\theta(u-v)\big)d\theta.
\end{equation*}
From Theorem \ref{pu}, we have the following Proposition.
\begin{proposition}\label{j1}
Under the assumptions given above, it holds
\begin{equation}\label{j11}
\|J_1(t,\cdot)\|_{q;\Omega}\leqslant C(\delta)(1+t)^{-\frac{2}{(m-1)(m+1)}+\frac{\delta}{q}}
\end{equation}
for any $\delta\in(0,1)$ and $q\in[1,\infty]$.
\end{proposition}
\begin{proof}
From \eqref{pro4}, \eqref{gest} and \eqref{wlr}, we have
\begin{equation*}
\|J_1\|_1=\int_\Omega|J_1|{\rm d}x\leqslant C(1+t)^{-\frac{2}{(m-1)(m+1)}}\int_\mathbb R\big((1-g)g+\partial_1g\big){\rm d}x_1\leqslant C(1+t)^{-\frac{2}{(m-1)(m+1)}+\delta}
\end{equation*}
for any $0<\delta<1$. On the other hand,
\begin{equation*}
\|J_1\|_\infty\leqslant C\big(\|\nabla u_l\|_\infty+\|\nabla u_r\|_\infty+\|w_l\|_\infty+\|w_r\|_\infty\big)\leqslant C(1+t)^{-\frac{2}{(m-1)(m+1)}}.
\end{equation*}
For $1<q<\infty$, \eqref{j11} follows from H\"{o}lder's inequality.
\end{proof}

\vskip 0.2in

Next, we will discuss the time-decay rate of $\phi$. Multiplying \eqref{phi} by $|\phi|^{r-2}\phi$ with $r\geqslant2$ implies
\begin{equation}\label{phi1}
\begin{aligned}
&\frac{1}{r}\frac{\partial}{\partial t}\big(|\phi|^r\big)\overbrace{-{\rm div}\big(|\nabla(\phi+\tilde u)|^{m-1}\nabla(\phi+\tilde u)-|\nabla\tilde u|^{m-1}\nabla\tilde u\big)|\phi|^{r-2}\phi}^{I_{11}}\\
&=\underbrace{-{\rm div}\big(f(\phi+\tilde u)-f(\tilde u)\big)|\phi|^{r-2}\phi}_{I_{12}}+\underbrace{\big({\rm div}(|\nabla\tilde u|^{m-1}\nabla\tilde u)-J_2\big)|\phi|^{r-2}\phi}_{I_{13}}-J_1|\phi|^{r-2}\phi.
\end{aligned}
\end{equation}
Note that
\begin{equation*}
I_{11}=\sum_{i=1}^N\partial_i(\cdots)+(r-1)\big(|\nabla(\phi+\tilde u)|^{m-1}\nabla(\phi+\tilde u)-|\nabla\tilde u|^{m-1}\nabla\tilde u\big)\cdot|\phi|^{r-2}\nabla\phi
\end{equation*}
with $(\cdots)=\big(|\nabla(\phi+\tilde u)|^{m-1}\partial_i(\phi+\tilde u)-|\nabla\tilde u|^{m-1}\partial_i\tilde u\big)|\phi|^{r-2}\phi$,
\begin{equation*}
\begin{aligned}
I_{12}&=\sum_{i=1}^N\partial_i(\cdots)-(r-1)\int_0^\phi\big(f_1^\prime(\eta+\tilde u)-f_1^\prime(\tilde u)\big)|\eta|^{r-2}\partial_1u^R{\rm d}\eta\\
&\qquad-(r-1)\int_0^\phi\big(f_1^\prime(\eta+\tilde u)-f_1^\prime(\tilde u)\big)|\eta|^{r-2}\partial_1(\tilde u-u^R){\rm d}\eta\\
&\qquad-(r-1)\sum_{i=2}^N\int_0^\phi\big(f_i^\prime(\eta+\tilde u)-f_i^\prime(\tilde u)\big)|\eta|^{r-2}\partial_i\tilde u{\rm d}\eta
\end{aligned}
\end{equation*}
with $(\cdots)=(r-1)\int_0^\phi\big(f_i(\eta+\tilde u)-f_i(\tilde u)\big)|\eta|^{r-2}{\rm d}\eta-(f(\phi+\tilde u)-f(\tilde u)\big)|\phi|^{r-2}\phi$,
\begin{equation*}
\begin{aligned}
I_{13}&=(1-g){\rm div}(|\nabla\tilde u|^{m-1}\nabla\tilde u-|\nabla u_l|^{m-1}\nabla u_l)|\phi|^{r-2}\phi\\
&\qquad+g{\rm div}(|\nabla\tilde u|^{m-1}\nabla\tilde u-|\nabla u_r|^{m-1}\nabla u_r)|\phi|^{r-2}\phi\\
&=\sum_{i=1}^N\partial_i(\cdots)+\partial_1g(|\nabla u_l|^{m-1}\partial_1u_l-|\nabla u_r|^{m-1}\partial_1u_r)|\phi|^{r-2}\phi\\
&\qquad-(r-1)\big(|\nabla\tilde u|^{m-1}\nabla\tilde u-(1-g)|\nabla u_l|^{m-1}\nabla u_l-g|\nabla u_r|^{m-1}\nabla u_r\big)\cdot|\phi|^{r-2}\nabla\phi
\end{aligned}
\end{equation*}
with $(\cdots)=\big(|\nabla\tilde u|^{m-1}\nabla\tilde u-(1-g)|\nabla u_l|^{m-1}\nabla u_l-g|\nabla u_r|^{m-1}\nabla u_r\big)|\phi|^{r-2}\phi$. Set
\begin{equation*}
\gamma=\frac{2}{(m-1)(m+1)}.
\end{equation*}
Since
\begin{equation*}
(r-1)\int_0^\phi\big(f_1^\prime(\eta+\tilde u)-f_1^\prime(\tilde u)\big)|\eta|^{r-2}\partial_1u^R{\rm d}\eta\geqslant c|\phi|^r\partial_1u^R\geqslant0
\end{equation*}
and
\begin{equation*}
\begin{aligned}
&\int_\Omega\Big|\partial_1g(|\nabla u_l|^{m-1}\partial_1u_l-|\nabla u_r|^{m-1}\partial_1u_r)|\phi|^{r-2}\phi\Big|{\rm d}x\\
&\hspace{1cm}\leqslant C(1+t)^{-1}\int_\Omega\big(|\nabla u_l|^m+|\nabla u_r|^m\big)|\phi|^{r-1}{\rm d}x\\
&\hspace{2cm}\leqslant C(1+t)^{-1}\big(\|u_l\|_{mr}^m+\|u_r\|_{mr}^m\big)\|\phi\|_r^{r-1}\leqslant C(1+t)^{-m\gamma-1}\|\phi\|_r^{r-1}
\end{aligned}
\end{equation*}
from \eqref{gest}, Theorem \ref{pu} and $rm>2$, we have, by integrating \eqref{phi1} with respect to $x$ over $\Omega$, that
\begin{equation}\label{phi2}
\begin{aligned}
&\frac{1}{r(r-1)}\frac{{\rm d}}{{\rm d}t}\|\phi\|_r^r+\overbrace{c\int_\Omega\big(|\nabla(\phi+\tilde u)|^{m-1}\nabla(\phi+\tilde u)-|\nabla\tilde u|^{m-1}\nabla\tilde u\big)\cdot|\phi|^{r-2}\nabla\phi{\rm d}x}^{I_{21}}\\
&+c\int_\Omega|\phi|^r\partial_1u^R{\rm d}x\leqslant C\big((1+t)^{-\gamma}\|\phi\|_r^r+(1+t)^{-m\gamma-1}\|\phi\|_r^{r-1}
+\|J_1\|_r\|\phi\|_r^{r-1}\big)\\
&+\underbrace{\int_\Omega\big(|\nabla\tilde u|^{m-1}\nabla\tilde u-(1-g)|\nabla u_l|^{m-1}\nabla u_l-g|\nabla u_r|^{m-1}\nabla u_r\big)\cdot|\phi|^{r-2}\nabla\phi{\rm d}x}_{I_{22}}.
\end{aligned}
\end{equation}
From \eqref{ab1}, we have
\begin{equation}\label{phi3}
I_{21}\geqslant c\int_\Omega\big(|\nabla(\phi+\tilde u)|^{m-1}+|\nabla\phi|^{m-1}+|\nabla\tilde u|^{m-1}\big)|\phi|^{r-2}|\nabla\phi|^2{\rm d}x.
\end{equation}
On the other hand, using \eqref{ab2} and \eqref{tildeu}, it holds
\begin{equation*}
\begin{aligned}
I_{22}&\leqslant \int_\Omega\Big((1-g)(|\nabla\tilde u|^{m-1}+|\nabla u_l|^{m-1})|\nabla(\tilde u-u_l)|\\
&\hspace{3cm}+g(|\nabla\tilde u|^{m-1}+|\nabla u_r|^{m-1})|\nabla(\tilde u-u_r)|\Big)|\phi|^{r-2}|\nabla\phi|{\rm d}x\\
&\leqslant \int_\Omega\big(|\nabla\tilde u|^{m-1}+|\nabla u_l|^{m-1}+|\nabla u_r|^{m-1}\big)\\
&\hspace{3cm}\big((1-g)g|\nabla(u_l-u_r)|+|\partial_1g|\big)|\phi|^{r-2}|\nabla\phi|{\rm d}x.
\end{aligned}
\end{equation*}
Then, it follows from Young's inequality, Theorem \ref{pu}, H\"{o}lder's inequality and \eqref{gest} that
\begin{equation}\label{phi4}
\begin{aligned}
I_{22}&\leqslant\varepsilon\int_\Omega\big(|\nabla\tilde u|^{m-1}+|\nabla\phi|^{m-1}\big)|\phi|^{r-2}|\nabla\phi|^2{\rm d}x\\
&\qquad+C\Big((1+t)^{-2\gamma-(m-1)\tilde\gamma}+(1+t)^{-\frac{2}{m-1}}\Big)\int_\Omega|\phi|^{r-2}(1-g)g{\rm d}x\\
&\qquad+C(1+t)^{-(m-1)\tilde\gamma}\int_\Omega|\phi|^{r-2}|\partial_1g|^2{\rm d}x+C(1+t)^{-\frac{2}{m}}\int_\Omega|\phi|^{r-2}|\partial_1g|^\frac{m+1}{m}{\rm d}x\\
&\leqslant\varepsilon\int_\Omega\big(|\nabla\tilde u|^{m-1}+|\nabla\phi|^{m-1}\big)|\phi|^{r-2}|\nabla\phi|^2{\rm d}x+C(1+t)^{-\varsigma}\|\phi\|_r^{r-2},
\end{aligned}
\end{equation}
where $\varepsilon>0$ is a small constant, $\tilde\gamma=\min\{1,\gamma\}=1$ and
\begin{equation*}
\begin{aligned}
\varsigma&=\min\{2\gamma+m-1-\frac{2\delta}{r},\quad\frac{2}{m-1}-\frac{2\delta}{r},\quad
m+1-\frac{2}{r},\quad\frac{m+3}{m}-\frac{2}{r}\}\\
&=m+1-\frac{2}{r}
\end{aligned}
\end{equation*}
for $\delta>0$ sufficiently small. Here, we used the fact that
\begin{equation*}
\|\nabla\tilde u\|_\infty\leqslant C(1+t)^{-\tilde\gamma}
\end{equation*}
from \eqref{tildeu}, \eqref{pro5} and \eqref{gest}. Thus, comparing (\ref{phi2}--\ref{phi4}), we have
\begin{equation}\label{phi5}
\begin{aligned}
&\frac{{\rm d}}{{\rm d}t}\|\phi\|_r^r+\int_\Omega\big(|\nabla(\phi+\tilde u)|^{m-1}+|\nabla\phi|^{m-1}+|\nabla\tilde u|^{m-1}\big)|\phi|^{r-2}|\nabla\phi|^2{\rm d}x\\
&+\int_\Omega|\phi|^r\partial_1u^R{\rm d}x\leqslant C\big((1+t)^{-\gamma}\|\phi\|_r^r
+(1+t)^{-\gamma+\delta}\|\phi\|_r^{r-1}+(1+t)^{-\varsigma}\|\phi\|_r^{r-2}\big).
\end{aligned}
\end{equation}
where we used \eqref{j11}. Note that $\varsigma>m>1$.

\vskip 0.2in

If we choose $r=2$, then \eqref{phi5} becomes
\begin{equation}\label{phi6}
\begin{aligned}
&\frac{{\rm d}}{{\rm d}t}\|\phi\|^2+\int_\Omega\big(|\nabla(\phi+\tilde u)|^{m-1}+|\nabla\phi|^{m-1}+|\nabla\tilde u|^{m-1}\big)|\nabla\phi|^2{\rm d}x+\int_\Omega|\phi|^2\partial_1u^R{\rm d}x\\
&\qquad\qquad\leqslant C\big((1+t)^{-2\gamma+\varsigma+2\delta}\|\phi\|^2+(1+t)^{-\varsigma}\big).
\end{aligned}
\end{equation}

\begin{lemma}\label{dt}
If a nonnegative function $y=y(t)$ satisfies $y(0)=0$ and
\begin{equation}\label{yt}
\frac{{\rm d}y}{{\rm d}t}\leqslant C_1(1+t)^{-\alpha}y+C_2(1+t)^{-\beta},
\end{equation}
where $C_1,C_2,\alpha>0$ and $\beta>1$. Then
\begin{equation}\label{yt2}
y(t)\leqslant\frac{C_2}{\beta-1}e^{\frac{C_1}{\alpha-1}}
\end{equation}
for any $t>0$.
\end{lemma}
\begin{proof}
Denote $C_0=\frac{C_1}{\alpha-1}$. Multiplying \eqref{yt} with $e^{C_0(1+t)^{-\alpha+1}}$, we have
\begin{equation*}
\frac{{\rm d}}{{\rm d}t}\left(e^{C_0(1+t)^{-\alpha+1}}y\right)\leqslant C_2(1+t)^{-\beta}e^{C_0(1+t)^{-\alpha+1}}\leqslant C_2(1+t)^{-\beta}e^{C_0}.
\end{equation*}
Integrating over $(0,t)$ implies
\begin{equation*}
y(t)\leqslant e^{-C_0(1+t)^{-\alpha+1}}\left(y(0)+\frac{C_2}{\beta-1}e^{C_0}\right)\leqslant\frac{C_2}{\beta-1}e^{C_0},
\end{equation*}
which is \eqref{yt2}
\end{proof}
Noting that $-2\gamma+\varsigma+2\delta<0$ and applying Lemma \ref{dt} to \eqref{phi6}, we can conclude
\begin{equation}\label{phi7}
\|\phi(t,\cdot)\|\leqslant C
\end{equation}
for any $t>0$.

\vskip 0.2in

For the case that $r>2$, we firstly suppose $r\leqslant m+3$. Multiplying \eqref{phi5} by $(1+t)^\beta$ for some constant $\beta>0$ and integrating the resultant equation over $(0,T)$ yield that
\begin{equation}\label{phi8}
(1+T)^\beta\|\phi(T,\cdot)\|_r^r+\int_0^T(1+t)^\beta\big\|\nabla|\phi|^\frac{r+m-1}{m+1}\big\|_{m+1}^{m+1}{\rm d}t\leqslant C\int_0^T(1+t)^{\beta-\nu}\|\phi\|_r^r{\rm d}t+C,
\end{equation}
where we used Young's inequality and $\nu\leqslant\min\{1,\gamma-\frac{\beta+1-\gamma}{r-1}-\delta, \varsigma-\frac{2(\beta+1-\varsigma)}{r-2}-\delta\}$ and $C$ is independent of $T$. Let $l_k,k=1,2,\cdots$ be a positive sequence to be determined below. Since $r\leqslant m+3$, it follows from Theorem 1.4 in \cite{hua21} that
\begin{equation}\label{phi9}
\begin{aligned}
&C\int_0^T(1+t)^{\beta-\nu}\|\phi\|_r^r{\rm d}t=C\int_0^T(1+t)^{\beta-\nu}\Big\||\phi|^\frac{r+m-1}{m+1}\Big\|_{r\frac{m+1}{r+m-1}}^{r\frac{m+1}{r+m-1}}{\rm d}t\\
&\leqslant C\sum_{k=0}^{N-1}\int_0^T(1+t)^{\beta-\nu}
\big\|\nabla|\phi|^\frac{r+m-1}{m+1}\big\|_{m+1}^{r\frac{m+1}{r+m-1}\theta_k}
\Big\||\phi|^\frac{r+m-1}{m+1}\Big\|_{l_1\frac{m+1}{r+m-1}}^{r\frac{m+1}{r+m-1}(1-\theta_k)}{\rm d}t\\
&=C\sum_{k=0}^{N-1}\int_0^T(1+t)^{\beta-\nu}
\big\|\nabla|\phi|^\frac{r+m-1}{m+1}\big\|_{m+1}^{r\frac{m+1}{r+m-1}\theta_k}\|\phi\|_{l_1}^{r(1-\theta_k)}{\rm d}t\\
&\leqslant \frac{1}{2}\int_0^T(1+t)^\beta\big\|\nabla|\phi|^\frac{r+m-1}{m+1}\big\|_{m+1}^{m+1}{\rm d}t
+C(1+T)^{\beta-\frac{(r+m-1)\nu+1}{r+m-1-r\theta_0}},
\end{aligned}
\end{equation}
where we used \eqref{phi7} for $l_1=2$ and
\begin{equation*}
\theta_k=\frac{(r-2)(r+m-1)}{r(r+m-1)-2r\frac{k-m}{k+1}}.
\end{equation*}
Here, since we only have the bound of $\|\phi\|$, we need $r\leqslant m+3$ so that $l_1\frac{m+1}{r+m-1}\geqslant1$. Comparing \eqref{phi8} and \eqref{phi9}, we have
\begin{equation*}
(1+T)^\beta\|\phi(T,\cdot)\|_r^r\leqslant C\left((1+T)^{\beta-\frac{r+3m-1}{3m+1}\nu+1}+1\right).
\end{equation*}
Since $1<m\leqslant\frac{3}{2}$, we can choose $\beta=\frac{r-2}{3m+1}$, then $\nu=1$ and
\begin{equation}\label{phi10}
\|\phi(T,\cdot)\|_r^r\leqslant C(1+T)^{-\frac{r-2}{3m+1}}
\end{equation}
for $2<r\leqslant m+3$.

Now we have $\|\phi(t,\cdot)\|_r^r\leqslant C$ for any $t>0$, so that we can suppose $r\leqslant m^2+3m+4$ by using $l_2=m+3$ instead of $l_1$ in \eqref{phi9} and further obtain \eqref{phi10} for $m+3<r\leqslant m^2+3m+4$. Repeating this progress, we can obtain \eqref{res3} and complete the proof of Theorem 3.

\vskip 0.2in

\begin{rem}
Since our $\phi$ does not vanish in the direction of $x_i,i=2,\cdots,N$, we can only use a Gagliado-Nirenburg (G-N) type inequality given in \cite{hua21}, instead of using the G-N inequality directly. Hence, the result in multi-dimension is not good as in 1-dimension.
\end{rem}

\section{\Large Appendix}
In this section, we will prove Lemma \ref{abp}. Let $a,b\in\mathbb R^N$ be two arbitrary vectors. Without loss of generality, we suppose $|a|\geqslant|b|$ in this section. Obviously, if $b=0$, Lemma \ref{abp} holds true, so we also suppose $|b|>0$.

\vskip 0.2in

To prove \eqref{ab2}, we denote
\begin{equation*}
\tilde a=\frac{a}{|b|},\quad \tilde b=\frac{b}{|b|},\quad \lambda=\frac{|a|}{|b|}.
\end{equation*}
Then $\lambda\geqslant1$, $|\tilde a|=\lambda|\tilde b|$ and
\begin{equation*}
(\lambda^{p-1}-1)|\tilde a|\leqslant(\lambda^{p-1}-\lambda^{p-2})|\tilde a|=|\lambda^p\tilde b-\lambda^{p-1}\tilde b|\leqslant|\lambda^{p-1}\tilde a-\lambda^{p-1}\tilde b|,
\end{equation*}
since $1<p<2$. Thus,
\begin{equation*}
|\lambda^{p-1}\tilde a-\tilde b|\leqslant|\tilde a-\tilde b|+(\lambda^{p-1}-1)|\tilde a|\leqslant(\lambda^{p-1}+1)|c-d|,
\end{equation*}
which immediately implies \eqref{ab2}.

\vskip 0.2in

Next, we will prove \eqref{ab1}. Obviously, \eqref{ab1} holds true for $a=\pm b$. We now suppose $a\neq\pm b$. It is easy to see that
\begin{equation*}
a\cdot b=|a||b|\cos\gamma,\qquad |a-b|^2=|a|^2+|b|^2-2|a||b|\cos\gamma,
\end{equation*}
where $\gamma$ is the included angle between $a$ and $b$. Set
\begin{equation*}
\alpha=\frac{|b|}{|a|},\qquad \beta=\cos\gamma,
\end{equation*}
then $\alpha\in(0,1],\beta\in[-1,1]$, $\alpha\neq1$ when $\beta=\pm1$.

We will firstly prove
\begin{equation}\label{bpp1}
\big(|a|^{q-1}a-|b|^{q-1}b\big)\cdot(a-b)\geqslant c(q)|a-b|^{q+1}.
\end{equation}
It is easy to see that
\begin{equation*}
\begin{aligned}
\frac{\big(|a|^{q-1}a-|b|^{q-1}b\big)\cdot(a-b)}{|a-b|^{q+1}}&=\frac{|a|^{q+1}-|a|^q|b|\cos\gamma-|a||b|^q\cos\gamma+ |b|^{q+1}}{\big(|a|^2+|b|^2-2|a||b|\cos\gamma\big)^{\frac{q+1}{2}}}\\
    &=\frac{1-\alpha\beta-\alpha^q\beta+\alpha^{q+1}}{\big(1+\alpha^2-2\alpha\beta\big)^{\frac{q+1}{2}}} :=f(\alpha,\beta).
\end{aligned}
\end{equation*}
In order to prove \eqref{ab1}, we only need to prove that there exists a positive constant $c$ such that $f(\alpha,\beta)\geqslant c$ for any $\alpha$ and $\beta$. It is easy to check that $f(\alpha,-1)=\frac{1+\alpha^q}{(1+\alpha)^q}>0,f(\alpha,1)=\frac{1-\alpha^q}{(1-\alpha)^q}>0$ for any $\alpha\in(0,1)$. Since $f$ is smooth, it remains to prove that $f(\alpha,\cdot)>0$ on extreme points in $(-1,1)$. Let $\partial_\beta f=0$, and direct calculation implies
\begin{equation*}
\beta=\frac{p\big(1+\alpha^{q+1}\big)-\alpha^{q-1}-\alpha^2}{(q-1)\alpha\big(1+\alpha^{q-1}\big)}.
\end{equation*}
Since $\beta<1$, we have
\begin{equation}\label{app1}
q<\alpha^{q-1}+(q-1)\alpha+(q-1)\alpha^q+\alpha^2-q\alpha^{q+1}=:h(\alpha).
\end{equation}
On the other hand, it holds
\begin{equation*}
h^\prime(\alpha)=(q-1)\alpha^{q-2}+q-1+q(q-1)\alpha^{q-1}+2\alpha-q(q+1)\alpha^q>0
\end{equation*}
for $\alpha\in(0,1)$, and $h(0)=0,h(1)=q$, which means $h(\alpha)\leqslant q$ for all $\alpha\in(0,1]$ and contradicts \eqref{app1}. This conclusion implies that $f(\alpha,\beta)$ is monotone with $\beta$, so that
\begin{equation*}
f(\alpha,\beta)\geqslant\min\{f(\alpha,-1),f(\alpha,1)\}\geqslant c(q)>0.
\end{equation*}
Then \eqref{bpp1} holds true.

To prove \eqref{ab1}, we only need
\begin{equation*}
|a-b|^{q-1}\geqslant c(q)(|a|^{q-1}+|b|^{q-1}).
\end{equation*}
Obviously,
\begin{equation*}
\frac{|a-b|^{q-1}}{|a|^{q-1}+|b|^{q-1}}=\frac{\big(1+\alpha^2-2\alpha\beta\big)^{\frac{q-1}{2}}}{1+\alpha^{q-1}}
:=g(\alpha,\beta).
\end{equation*}
Since $g(\alpha,-1)=\frac{(1+\alpha)^{q-1}}{1+\alpha^{q-1}}>0$, $g(\alpha,1)=\frac{(1-\alpha)^{q-1}}{1+\alpha^{q-1}}>0$, and
\begin{equation*}
\partial_\beta g=-\frac{(q-1)\alpha}{1+\alpha^{q-1}}\big(1+\alpha^2-2\alpha\beta\big)^{\frac{q-3}{2}}<0,
\end{equation*}
we have
\begin{equation*}
g(\alpha,\beta)\geqslant\min\{g(\alpha,-1),g(\alpha,1)\}\geqslant c(q)>0,
\end{equation*}
which completes the proof.

\end{document}